\pdfoutput=1 

\documentclass[12pt,twoside]{preprint}

\title{On the complex structure of symplectic quotients}
\author{Xiangsheng Wang}
\currentaddress{School of Mathematics, Shandong University, Jinan, Shandong {\rm 250100}, China}
\institute{BICMR, Peking University, Beijing {\rm 100871}, China\hspace{3.5em} \email{wangxs1989@gmail.com}}
\date{\today}

\usepackage[english]{babel}
\usepackage[stretch=10,final]{microtype}
\usepackage[a4paper,left=2.5cm,right=2.5cm,top=4cm,bottom=2.5cm,marginparwidth=2cm]{geometry}
\usepackage{color}
\usepackage{titlesec}
\usepackage{enumitem}
\usepackage[olditem,oldenum]{paralist}
\usepackage{perpage}
\usepackage[super]{nth}
\usepackage{todonotes}
\usepackage{fancyhdr}
\usepackage[title]{appendix}
\usepackage[colorlinks,linkcolor=blue,bookmarksnumbered,pagebackref,final]{hyperref}
\usepackage[]{lipsum}
\usepackage[alphabetic,initials,msc-links]{amsrefs}

\usepackage{mathtools}
\usepackage{amsthm}
\usepackage{amssymb}
\usepackage{tikz-cd}

\usepackage[osf,theoremfont]{newtxtext} \usepackage[upint,smallerops]{newtxmath}
\usepackage[cal=euler,scr=rsfso,frak=euler]{mathalfa}
\usepackage[T1]{fontenc}

\usepackage[normalem]{ulem}

\makeatletter
\AtBeginDocument{
  \hypersetup{
    pdftitle = {\@title},
    pdfauthor = {\@author}
  }
}
\makeatother

\MakePerPage{footnote}

\DeclareSymbolFont{sfoperators}{OT1}{cmss}{m}{n}
\DeclareSymbolFontAlphabet{\mathsf}{sfoperators}
\makeatletter
\renewcommand{\operator@font}{\mathgroup\symsfoperators}
\makeatother

\numberwithin{equation}{section}
\swapnumbers                            \newtheoremstyle{plain}
{2ex plus 1ex minus .2ex}   {\medskipamount}   {\slshape}  {}       {\bfseries\tlfstyle} {.}         {5pt plus 1pt minus 1pt} {}          

\newtheorem{theorem}[subsection]{Theorem}
\newtheorem{lemma}[subsection]{Lemma}

\newtheorem{proposition}[subsection]{Proposition}
\newtheorem*{claim*}{Claim}
\newtheorem*{lemma*}{Lemma}

\newtheoremstyle{definition}
{2ex plus 1ex minus .2ex}   {\medskipamount}   {}  {}       {\bfseries\tlfstyle} {.}         {5pt plus 1pt minus 1pt} {}          \theoremstyle{definition}

\newtheorem{remark}[subsection]{Remark}
\newtheorem*{remark*}{Remark}
\newtheorem*{assumption*}{Assumption}

\makeatletter
\let\newtitle\@title
\let\newauthor\@author
\let\newdate\@date
\makeatother

\titleformat{\section}{\normalfont\large\bfseries\tlfstyle}{\thesection}{0.6em}{}
\titleformat{\subsection}[runin]{\normalfont\normalsize\bfseries\tlfstyle}{\thesubsection}{0.4em}{}

  \newcommand{\kk}{\mathfrak{k}}
\newcommand{\kg}{\mathfrak{g}}
\newcommand{\kt}{\mathfrak{t}}

\newcommand{\kz}{\mathfrak{z}}
\newcommand{\kr}{\mathfrak{r}}
\newcommand{\cO}{\mathcal{O}}
\newcommand{\cL}{\mathcal{L}}
\newcommand{\cR}{\mathcal{R}}
\newcommand{\cF}{\mathcal{F}}
\newcommand{\sL}{\mathscr{L}}
\newcommand{\sX}{\mathscr{X}}
\newcommand{\rT}{\mathrm{T}}
\newcommand{\rM}{\mathrm{M}}
\newcommand{\rs}{\mathrm{s}}
\newcommand{\opi}{\mathrm{i}}
\newcommand{\rss}{\mathrm{ss}}
\newcommand{\bZ}{\mathbb{Z}}
\newcommand{\bC}{\mathbb{C}}
\newcommand{\bR}{\mathbb{R}}
\newcommand{\bP}{\mathbb{P}}
\newcommand{\bQ}{\mathbb{Q}}
\newcommand{\bfv}{\mathbf{v}}
\newcommand{\bfw}{\mathbf{w}}

\newcommand{\bfT}{\mathbf{T}}
\newcommand{\kah}{K\"ahler}
\newcommand{\imp}{\mathrm{impl}}
\newcommand{\ms}{\mathrm{ms}}
\newcommand{\ums}{G_{\ms}}
\newcommand{\uip}{(\mathrm{T}^*K)_{\imp}}
\newcommand{\prin}{\mathrm{prin}}
\newcommand{\dsl}{\mathbin{/\mkern-5mu/}}

\newcommand{\KG}{K\backslash G}
\newcommand{\bKG}{\partial_{\infty} (K\backslash G)}
\newcommand{\norm}[1]{\left\Vert#1\right\Vert}

\renewcommand{\Im}{\operatorname{Im}\mathop{}\!}
\renewcommand{\emptyset}{\varnothing}

\DeclareMathOperator{\diff}{d}

\DeclareMathOperator{\cst}{J}
\DeclareMathOperator{\prj}{P}
\DeclareMathOperator{\opz}{z}
\DeclareMathOperator{\opsh}{sh}
\DeclareMathOperator{\opc}{c}
\DeclareMathOperator{\id}{id}
\DeclareMathOperator{\inc}{in}
\DeclareMathOperator{\Inc}{IN}
\DeclareMathOperator{\opp}{p}
\DeclareMathOperator{\opq}{q}
\DeclareMathOperator{\conv}{conv}

\begin{document}

\maketitle
\begin{abstract}
  Let $K$ be a compact group. For a symplectic quotient $M_{\lambda}$ of a compact
  Hamiltonian \kah\ $K$-manifold, we show that the induced complex structure on
  $M_{\lambda}$ is locally invariant when the parameter $\lambda$ varies in
  $\mathrm{Lie}(K)^*$. To prove such a result, we take two different
  approaches:
  \begin{inparaenum}[(i)]
  \item by using the complex geometry properties of the symplectic implosion
    construction;
  \item by investigating the variation of GIT quotients.
  \end{inparaenum}
\end{abstract}

\section{Introduction}

\subsection{}
For symplectic manifolds with Hamiltonian action, the symplectic reduction has
always been a powerful tool to investigate properties on such manifolds since
Marsden and Weinstein~\cite{Marsden:1974aa} introduced it several decades ago.
By using this technique, we can obtain a new symplectic manifold called
reduction manifold or symplectic quotient. The symplectic reduction establishes
a natural correspondence between the symplectic quotient and the original
manifold. As a consequence, if some quantities can be defined on both manifolds,
people can often find some interesting relations by comparing them. As an
illustration of this general idea, we would like recall two classical examples.

\texttt{Riemann-Roch numbers.} Since each symplectic manifold can also be viewed
as an almost complex manifold, we can define Riemann-Roch (RR) numbers for the
original manifold as well as the symplectic quotient. If we try to compare these
two RR numbers, we will encounter the famous geometric quantization conjecture,
c.f.~\cite{Guillemin:1982vb}. Because RR numbers can contain some representation
theoretic information, the merit of such a comparison is that the representation
theoretic information of the original manifold can be recovered from that of
the symplectic quotient.

\texttt{\kah\ metrics.} If the manifold is a \kah\ manifold, the symplectic
quotient can inherit a \kah\ structure automatically. In this case, \kah\
metrics are natural quantities that can be utilized for comparison. To obtain
specific results, we usually need some assumptions on the metric of the origin
manifold. For example, in~\cite[Ch. \S 7]{Futaki:1988aa}, Futaki obtains a
formula for the metric on the symplectic quotient supposing that the metric on
the original manifold has positive Ricci curvature. More recent development
along the same line appears in~\cite{La-Nave:2016en}. In this paper, La Nave and
Tian use a special equation to restrict the metric on the original \kah\
manifold. Then the metrics on symplectic quotients must satisfy the famous
\kah-Ricci flow equation. Moreover, in such a situation, they discover that the
symplectic reduction construction can be reversed in a certain sense, that is,
the solution to the \kah-Ricci flow equation can be used to construct the
special metric on the original manifold.

\subsection{}
The symplectic reduction has a nice feature, that is, we can usually obtain a
family of symplectic quotients rather than a single one. Acute readers may have
noticed that the above two examples have used this family of symplectic
quotients. As observed from these two examples, if we want to recover
information on the original manifold from the symplectic reduction, only input
from a single symplectic quotient is insufficient, in other words, we should
study the family of symplectic quotients. For this purpose, it is beneficial to
find some methods to compare the different symplectic quotients in the family.
In this paper, we are going to deal with one such comparison problem. For a
precise statement, we need some notations.

Let $K$ be compact group and $\kk$ be its Lie algebra. Choose a maximal torus
$T\subseteq K$ (a Cartan subalgebra $\kt\subset\kk$ resp.) and a closed positive
Weyl chamber $\kt^*_+ \subseteq \kt^*$. By the root spaces decomposition,
$\kt^*$ is naturally included in $\kk^*$. Let $M$ be a \kah\ manifold with a
holomorphic Hamiltonian $K$-action. We denote the symplectic quotient at
$\lambda \in \kk^*$ by $M_{\lambda}=M\dsl_{\lambda}K$, where $\kk^*$ is the dual
of the Lie algebra of $K$. It is well known that one can also give a complex
structure on $M\dsl_{\lambda}K$. In the following result, we compare the
symplectic quotients by their complex structures. For the sake of simplicity, in
this paper, we always assume that if $\lambda$ is a regular value of the moment
map $\mu$ of $M$, $K$ acts on $\mu^{-1}(\lambda)$ freely. In other words,
$M_{\lambda}$ is smooth.

\begin{theorem}
  \label{thm:intro1}
  Let $M$ be a compact Hamiltonian \kah\ $K$-manifold. Suppose that $\lambda \in \kt^*_+\subseteq
  \kk^*$ is a regular value of the moment map. Let $\sigma$ be the face of the Weyl
  chamber where $\lambda$ lies in. Then for any $\lambda'\in \sigma$ lying in a sufficiently small
  neighborhood of $\lambda$, $M_{\lambda'}$ is biholomorphic to $M_{\lambda}$.
\end{theorem}
Roughly speaking, the above result asserts that the complex structures on
$M_{\lambda}$ do not change when $\lambda$ varies in a small neighborhood. In the following,
we also refer to this property of symplectic quotients as the \emph{local
  invariance of complex structures} (with respect to the parameter $\lambda$). For
abelian group actions, it is known that symplectic quotients satisfy such a
property. But for the non-abelian case, it seems that the result similar to the
above theorem is not documented in literature explicitly as far as we concerned.
One aim of this paper is to fill this gap.

As we have explained at the beginning of this subsection, we will treat the
above theorem as a premise to study the relation between some quantities defined
both on the original manifold and symplectic quotients. In this specific
situation, one possible direction along this line is trying to generalize La
Nave-Tian's correspondence of \kah\ metrics mentioned before to non-abelian
group actions based on the above result.

As to the proof of Theorem~\ref{thm:intro1}, we will take two different
approaches: one is of more symplectic geometry flavor while another is of more
algebraic geometry flavor. Although looked different ostensibly, we find that
both of them are linked to a simple idea: reducing a non-abelian problem to an
abelian one.

\subsection{}
For the symplectic geometry method for Theorem~\ref{thm:intro1}, we use a
construction called \emph{symplectic implosion} due to Guillemin, Jeffery and
Sjamaar~\cite{Guillemin:2002sx}. Roughly speaking, such a construction enables
us to construct a new symplectic space with an abelian group action to
substitute the original one. We can hint the usage of this construction by the
following result in~\cite{Guillemin:2002sx}, c.f.\ Theorem~\ref{thm:symquo}: the
symplectic quotient of a manifold, which is generally obtained from a
non-abelian reduction, is symplectically isomorphic to the abelian symplectic
quotient of its symplectic implosion.

To prove Theorem~\ref{thm:intro1}, we need to find a \kah\ manifolds refinement
of the result of Guillemin et al.\ mentioned above. More precisely, we have the
following result which compares the complex structures (in fact, the \kah\
structures) on a symplectic quotient inherited from the original manifold and
its symplectic implosion respectively.

\begin{theorem}
  \label{thm:introeq}
  Let $M$ be a Hamiltonian \kah\ $K$-manifold and let $M_{\imp}$ be the
  symplectic implosion of $M$. For any $\lambda \in \kt^*_+\subseteq \kk^*$
  being a regular value of the moment map of $M$, the two symplectic quotients,
  $M\dsl_{\lambda}K$ and $M_{\imp}\dsl_{\lambda}T$, are naturally isomorphic as
  \kah\ manifolds.
\end{theorem}

As a side comment, we notice that, in~\cite[Theorem~3.8]{Safronov:2017gd},
Safronov proves a result similar to Theorem~\ref{thm:introeq} for algebraic
varieties with a complex symplectic structure. Compared to the method therein,
i.e.\ derived symplectic geometry, we follow a rather elementary calculation
approach for the proof of Theorem~\ref{thm:introeq}.

With Theorem~\ref{thm:introeq}, we can complete the first proof of
Theorem~\ref{thm:intro1} by using the local invariance of complex structures
for abelian symplectic quotients.

\subsection{}
Symplectic quotients have an intimate relation with geometric invariant theory
(GIT) quotients, c.f.~\cites{Kirwan:1984st,Ness:1984rr,Sjamaar:1995aa} for
smooth manifolds and~\cites{Heinzner:1994bn,Heinzner:1999aa} for general complex
spaces. Considering this, we would like to give some GIT explanations
of Theorem~\ref{thm:intro1} before turning to the second proof of it. More
concretely, we use the symplectic implosion again to show the following result,
referred to Proposition~\ref{prop:ss} for a more precise statement and a
discussion about the stability condition used here.

\begin{proposition}
  \label{prop:sseq}
  Let $\lambda,\lambda'$ be in the same positions as in Theorem \ref{thm:intro1}. Denote
  the coadjoint orbit going through $\lambda$ ($\lambda'$) by $\cO_{\lambda}$ ($\cO_{\lambda'}$). Then
  the semi-stable point sets of $M \times \cO_{\lambda}^*$\footnote{The asterisk here
    means that we use the negative canonical symplectic form on a coadjoint
    orbit.} and $M\times \cO_{\lambda'}^*$ coincide.
\end{proposition}

\subsection{}
Our second approach to Theorem~\ref{thm:intro1} is based on the Kirwan-Ness
theorem, namely, the coincidence of GIT quotients and symplectic quotients. Due
to this theorem, if we are able to show Proposition~\ref{prop:sseq} without using
symplectic implosion, we can find another proof of Theorem~\ref{thm:intro1}.
Therefore, what we have actually done is using some complex algebraic geometry
techniques to reprove Proposition~\ref{prop:sseq}.

The main ideas and tools behind this method come from the variation of GIT
quotients theory (vGIT) as developed by~\cite{Dolgachev:1998wv} and
\cite{Thaddeus:1996tt}. In the analytic language, vGIT discusses the variation
of the symplectic quotients when the symplectic form (or the moment map) on
the original manifold changes. From this viewpoint, it is quite reasonable to
expect that Proposition~\ref{prop:sseq} will follow from general vGIT results.
In fact, Proposition~\ref{prop:sseq} can be seen as a local and weaker version
of the following result, c.f. Theorem~\ref{thm:dol}~\&~\ref{thm:finite}.
\begin{theorem}
  \label{thm:global-prop}
  Let $K_{\sigma} \subseteq K$ be the isotropic subgroup associated to face $\sigma$. For any
  $\lambda\in \sigma$, let $(M \times K/K_{\sigma})^{\rss}_{\lambda} \subseteq M \times K/K_{\sigma}$ be the
  semi-stable point set associated to $\lambda$. One concludes that there is a finite
  partition of $\sigma = \sqcup^N_{i=1} \sigma_i$, such that for any $\sigma_i$ which is not a
  ``wall'', and $\lambda,\lambda'$ lying in the interior of $\sigma_i$, $(M \times
  K/K_{\sigma})^{\rss}_{\lambda} = (M \times K/K_{\sigma})^{\rss}_{\lambda'}$.
\end{theorem}

In fact, in the above theorem, for the case $\sigma = \kt^{\circ}_+$, by using the Sard's
theorem, we can see that $\sigma_i$ must not be a wall. To show the theorem, we
need to generalize some vGIT results in \cite{Dolgachev:1998wv} and
\cite{Thaddeus:1996tt} from the projective algebraic varieties to the general
\kah\ manifolds. Especially, following~\cite{Dolgachev:1998wv}, we discuss the
relation between the stability condition and a special numerical function
$\rM^{\bullet}(x)$. The properties of $\rM^{\bullet}(x)$ play a major role in the proof of the
above theorem.

\subsection{}
In the final part of this paper, we discuss an interesting relation between
these two approaches of Theorem~\ref{thm:intro1}. Besides the similar philosophy
behind them, one can even find a precise correspondence in a certain way. We observe that
to reduce the vGIT problem to the abelian group action case,
in~\cite{Thaddeus:1996tt}, Thaddeus uses a construction very similar to the
symplectic implosion. In fact, in the setting we concerned in this paper, i.e.
on the manifold $M \times \cO^*_{\lambda}$, the output of Thaddeus's construction is a
partial desingularization of the symplectic implosion, c.f.
Proposition~\ref{prop:ms-imp}~\&~\ref{prop:master-space}.

\subsection{}
This paper is organized as follows. In Section~\ref{sec:si}, we review some
materials about the symplectic implosion used in this paper. After that, we
prove Theorem~\ref{thm:introeq} and give the first proof of
Theorem~\ref{thm:intro1} in Section~\ref{sec:1st-app}, using the symplectic
implosion. Section~\ref{sec:ss} is a transitional section, where we recall the
stability condition used in this paper and prove Proposition~\ref{prop:sseq}. In
Section~\ref{sec:vgit}, we discuss the properties of numerical function
$\rM^{\bullet}(x)$ and prove two vGIT type results
Theorem~\ref{thm:dol}~\&~\ref{thm:finite}. As a corollary, we obtain another
proof of Theorem~\ref{thm:intro1}. The last section,
Section~\ref{sec:comp-betw-two}, is devoted to the comparison of these two
approaches.

\subsection*{Acknowledgement.}
\label{sec:ack}
This work is supported by China Postdoctoral Science Foundation (Grant No.\
BX201700008). The author appreciates Prof. Gang Tian bringing the
paper~\cite{La-Nave:2016en} to his notice and encouraging him to consider its
non-abelian generalization, which is the major motivation for this paper. The
author also thanks the referees for very careful proofreading and helpful
comments, especially pointing out the important paper~\cite{Greb_2010ka} to him.

\section{Symplectic implosion}
\label{sec:si}

\subsection{}
In this section, we will review some backgrounds about the symplectic implosion.
Basically, all materials in this section follow closely with
\cite{Guillemin:2002sx}. Along the way, we also set up the assumptions and
notations used in the whole paper. We start with the symplectic geometry
features of the symplectic implosion, then turn to its complex geometry
properties.

\subsection{Symplectic aspects.}
Let $(M,\omega)$ be a connected symplectic manifold with a Hamiltonian group action
of a compact group $K$. Recall that an action is called \emph{Hamiltonian} if
there is a moment map $\mu$ for this action, which by definition is an
equivariantly map from $M$ to $\kk^*$, the dual of the Lie algebra of $K$, and
satisfies the following equation,
\begin{equation}
  \label{eq:mp}
  \diff \langle \mu, X \rangle = \iota_{X^M}\omega,
\end{equation}
where $X\in \kk$ and $X^M$ is the vector field induced by $X$ on $M$ by the
infinitesimal group action. We should remind readers that sometimes the moment
map is defined to be $-\mu$ in literature.

In this paper, we fix once and for all a maximal torus $T$ in $K$, denoting the
corresponding Cartan subalgebra by $\kt$. Also, we fix a closed positive Weyl
chamber $\kt^*_+$ in the dual of $\kt$. By using the root spaces decomposition
of $\kk$, $\kt^*$ is identified as a subspace of $\kk^*$. Choose a face $\sigma \subseteq
\kt^*_+$\footnote{$\sigma$ is also called a wall in literature. But in this paper,
  the terminology ``wall'' is reserved for another concept used latter.}, all
points lying on $\sigma$ have the same isotropic group $K_{\sigma}$ under the coadjoint
action $\mathrm{Ad}^*(K)$ on $\kt^*$. An equivalence relation $\sim$ is introduced
for points in $\mu^{-1} (\sigma)$ as follows: for $x,y$ in $\mu^{-1}(\sigma)$, $x \sim y$ if and
only if $x = ky$ for some $k \in [K_{\sigma}, K_{\sigma}]$. The \emph{imploded
  cross-section}\footnote{This is the terminology used by Guillemin \textit{et
    al.},~\cite{Guillemin:2002sx}. In this paper, imploded cross-section and \emph{symplectic implosion}
  will be used interchangeably.} of $M$,
\cite[Definition~2.1]{Guillemin:2002sx}, is defined to be the quotient space:
$M_{\imp}\coloneq \mu^{-1}(\kt^*_+)/\sim$, with the quotient map denoted by $\pi:
\mu^{-1}(\kt^*_+) \rightarrow M_{\imp}$. Set-theoretically, the imploded cross-section can
be written as the following disjoint union,
\begin{equation}
  \label{eq:1}
  M_{\imp} = \coprod_{\sigma \in \Sigma} \mu^{-1}(\sigma) / [K_{\sigma},K_{\sigma}],
\end{equation}
where $\Sigma$ denotes the index set of the faces of $\kt^*_+$. Note that on $\Sigma$,
there is a natural partial order: $\sigma \le \tau $ if and only if $\sigma \subseteq \bar{\tau}$.

\begin{remark}
  \label{rk:sidef}
  About the imploded cross-section $M_{\imp}$, the following properties hold.
  \begin{enumerate}[label=\Alph*.]
  \item In general, $M_{\imp}$ is not a smooth manifold, but only can be a
    stratified symplectic space in the weak sense of \cite{Sjamaar:1991aa}. The
    quotient map $\pi$ is always proper. If we assume that $M$ is compact, so is
    $M_{\imp}$.
  \item By~\cite[Corollary~2.7]{Guillemin:2002sx}, every component appearing in
    the decomposition of (\ref{eq:1}) is a symplectic quotient. More precisely,
    one has
    \begin{equation}
      \label{eq:2}
      M_{\imp} = \coprod_{\sigma \in \Sigma} \mu^{-1}(\sigma) / [K_{\sigma},K_{\sigma}]= \coprod_{\sigma\in \Sigma}{\mu^{-1}\big(K_{\sigma}(\cup_{\tau \ge \sigma}\tau)\big)
      } \dsl [K_{\sigma},K_{\sigma}].\footnote{When taking the
        symplectic quotient with respect to $0$, the subscript of $\dsl$ is omitted.}
    \end{equation}
    Moreover, by the symplectic cross-section
    theorem,~\cite[Theorem~2.5]{Guillemin:2002sx}, for any $\sigma \in \Sigma$,
    $\mu^{-1}\big(K_{\sigma}(\cup_{\tau \ge \sigma}\tau)\big)$ is a smooth submanifold of $M$, which
    implies the singularity of $M_{\imp}$ is not too bad.
  \item The minimal face $\sigma$ satisfying $\mu(M)\subseteq \bar{\sigma}$ is called the
    \emph{principal face} for $M$, denoted by $\sigma_{\prin}$. The group action
    $[K_{\sigma_{\prin}}, K_{\sigma_{\prin}}]$ on $\mu^{-1}(\sigma_{\prin})$ is trivial actually,
    which means that the stratum
    \begin{equation*}
      \mu^{-1}(\sigma_{\prin})/[K_{\sigma_{\prin}},K_{\sigma_{\prin}}]=\mu^{-1}(\sigma_{\prin})
    \end{equation*}
    (called the \emph{principal cross-section} of $M$) in decomposition
    (\ref{eq:1}) must be smooth. In many cases, the principal face of $M$ is the
    interior of the positive Weyl chamber $(\kt^*_+)^{\circ}$. Especially, since
    in this paper we are mainly concerned the case that the regular value set of
    $\mu$ is non-empty, \uline{we will always assume that
      $\sigma_{\prin}=(\kt^*_+)^{\circ}$}.
  \end{enumerate}
\end{remark}

\subsection{}
By definition, $M_{\imp}$ inherits a $T$-action from the $K$-action of $M$.
Besides, the moment map $\mu$ on $M$ also induces a continuous map
$\mu_{\imp}:M_{\imp}\rightarrow \kt^*_+$ on the symplectic implosion. Although $M_{\imp}$ is
not a manifold in general, $(M_{\imp},\omega_{\imp},\mu_{\imp})$ can be seen as a
Hamiltonian $T$-space. This means that when restricted to a smooth
stratum, $\mu_{\imp}$ is just a moment map for the smooth $T$-action in the usual
sense. The symplectic reduction construction can also be generalized to this
case. The following theorem asserts that the ``$K$-reduction'' on $M$ is equal
to the ``$T$-reduction'' on $M_{\imp}$, which explains why the symplectic
implosion $(M_{\imp},\omega_{\imp},\mu_{\imp})$ can be seen as the abelianization of
$(M,\omega,\mu)$ as we have recalled in Introduction.
\begin{theorem}[{\cite[Theorem~3.4]{Guillemin:2002sx}}]
  \label{thm:symquo}
  For every $\lambda$ in a face $\sigma$ of $\kt^*_+$, the canonical map from $\mu^{-1}(\lambda)$
  to $\mu^{-1}_{\imp}(\lambda) = \mu^{-1}(\lambda)/[K_{\sigma},K_{\sigma}]$ induces a symplectic
  isomorphism between the two symplectic quotients $M \dsl_{\lambda} K \simeq M_{\imp}
  \dsl_{\lambda} T$.
\end{theorem}

\subsection{}
As an example of the above construction, we take a closer look at $\uip$, the
symplectic implosion of the cotangent bundle of $K$. Quite interesting, as we
will see, this special symplectic implosion has certain universal properties.

For the later usage, here we also set up the convention of the group actions on
$K$ and $\rT^*K$. The left and right action of $K$ on itself are denoted by
$\cL_g k \coloneq gk$, $\cR_g k \coloneq kg^{-1}$, where $k$, $g$ are elements of $K$.
We identity $\rT^*K$ ($\rT K$ resp.) with $K \times \kk^*$ ($K \times \kk$ resp.) by the
left translation to the identity element. Under such identification, the
Liouville form $\beta$ on $\rT^*K$ can be written as follows,
\begin{equation}
  \label{eq:kliou-form}
  \beta_{(k,\lambda)}(X,\xi) = \langle\lambda,X\rangle,\quad \text{for }(k,\lambda)\in K \times \kk^* \simeq \rT^*K,\, (X,\xi)\in \kk \times \kk^* \simeq \rT_{(k,\lambda)}(K \times \kk^*).
\end{equation}
Let $\omega = \diff \beta$ be the usual symplectic form on $\rT^*K$. Note that both $\cL$
and $\cR$ can be lifted to actions on $\rT^*K$. We denote the lifted actions by
the same notations: $\cL_g(k,\lambda)=(gk,\lambda)$, $\cR_g(k,\lambda)=(kg^{-1},g\lambda)$, where the
group action on $\kk^*$ is the coadjoint action $\mathrm{Ad}^*$. Both of $\cL$
and $\cR$ are Hamiltonian with respect to $\omega$, by (\ref{eq:kliou-form}), whose
moment maps are,
\begin{equation}
  \label{eq:kmm}
  \Phi_{\cL}(k,\lambda) = -k\lambda,\quad \Phi_{\cR}(k,\lambda) = \lambda.
\end{equation}

One can carry out the symplectic implosion construction for $\rT^*K$ using
either $\cL$ or $\cR$. Following the convention of~\cite{Guillemin:2002sx}, in
this paper, the imploded cross-section of $\rT^*K$, $\uip$, is always
constructed out of $\cR$ unless otherwise declared. The decomposition of
(\ref{eq:2}) here has a more explicit form,
\begin{equation}
  \label{eq:tkdc}
  \uip = \coprod_{\sigma \in \Sigma}\big(K \times (K_{\sigma} \cup_{\tau \ge \sigma}\tau)\big)\dsl [K_{\sigma},K_{\sigma}] = \coprod_{\sigma \in \Sigma} \frac{K}{[K_{\sigma},K_{\sigma}]} \times \sigma.
\end{equation}
Note that the principal cross-section here is $K \times (\kt^*_+)^{\circ}$ and all
components in the decomposition (\ref{eq:tkdc}) are smooth submanifolds. As the
symplectic implosion of $\rT^*K$, $\uip$ inherits a $K$-action from $\cL$ and a
$T$-action from $\cR$. By (\ref{eq:kmm}), the moment map for the $T$-action on
$\uip$ is $\Phi_{\imp,\cR}([k],\lambda) = \lambda$ for $[k]\in {K}/{[K_{\sigma},K_{\sigma}]}$. Moreover,
it's easy to check that the induced $K$-action is also Hamiltonian whose moment
map is $\Phi_{\imp,\cL}([k],\lambda) = -k\lambda$. By Theorem~\ref{thm:symquo}, one can
calculate the symplectic quotient of $\uip$ with respect to the $T$ action
explicitly. Recall that a well-known property of $\rT^*K$ is that its symplectic
quotients are coadjoint orbits,~\cite[\S~4.2]{Marsden:1974aa}. Therefore, the
symplectic quotient $\uip \dsl_{\lambda} T$ as a Hamiltonian $K$-manifold, is
naturally isomorphic to the coadjoint orbit $\cO^*_{\lambda}$.

A special feature of $\uip$ is the following universal property.
By~\cite[Theorem~4.9]{Guillemin:2002sx}, for any Hamiltonian $K$-manifold
$(M,\omega,\mu)$, $M_{\imp}$ as a Hamiltonian $K$-space can be constructed as follows,
\begin{equation}
  \label{eq:up}
  M_{\imp} = (M \times \uip) \dsl K,
\end{equation}
where the $K$-action on the product manifold is the diagonal action. For this
reason, $\uip$ is called the \emph{universal imploded cross-section}.

\subsection{}
\label{sub:group-ass}
Some comments about the assumption on the groups used in this paper. As we have
said, we always assume that $K$ is a compact group. Let $K'$ ba finite cover of
$K$. The $K$-action on $M$ also induces a Hamiltonian $K'$-action. And
symplectic quotients obtained with respect to $K$-action and $K'$-action are the
same. Since in this paper, we only interested in the properties of symplectic
quotients, from now on, we will assume that $K$ is isomorphic to a product of a
torus and a semi-simple simply connected group.

Moreover, in the remaining subsections of this section discussing the
complex geometry properties of a symplectic implosion, we will further assume
that the group $K$ is a semi-simple simply connected group to omit some
technicalities involved in the general case~\cite[the bottom of
p.~174]{Guillemin:2002sx}. We remark that such an assumption causes no loss of
generalities. Since, by~\cite[Lemma~2.4]{Guillemin:2002sx}, the symplectic
implosions of $M$ with respect to $K$ and $[K,K]$ are equal, one can use the
complex structure of $M_{\imp,[K,K]}$ to define the complex structure of
$M_{\imp,K}$.

\subsection{Complex aspects.}
\label{sub:complex-aspects}
In the rest of this paper, we assume that the Hamiltonian $K$-manifold
$(M,\omega,\mu)$ is endowed with a compatible integrable $K$-invariant complex
structure $\cst$, which means that $\cst$ preserves the symplectic structure and
$g(-,-) \coloneq \omega(-,\cst -)$ is a Riemannian metric on $M$. In other words,
$(M,g,\cst)$ is a \kah\ manifold and $\omega$ is the \kah\ form. Unlike the
symplectic structure on $M_{\imp}$, the \kah\ structure on $M_{\imp}$ does not
inherit from $M$ directly. To see this, consider the following fact: although
the principal cross-section $\mu^{-1}(\sigma_{\prin})$ of $M_{\imp}$ can be seen as a
smooth submanifold of $M$, in general, $\mu^{-1}(\sigma_{\prin})$ is not a complex
submanifold of $M$. To overcome this problem, one first defines the complex
structure on $\uip$ and the complex structure of $M_{\imp}$ is defined by using
the equality (\ref{eq:up}).

Following \cite[\S~6]{Guillemin:2002sx}, to define the complex structure on the
universal imploded cross-section $\uip$, we need to embed it into a
$K$-representation space as an affine subvariety. Let $\Lambda= \ker(\exp|_{\kt})$ be
the exponential lattice in $\kt$ and $\Lambda^*= \mathrm{Hom}_{\bZ}(\Lambda, \bZ)$ be the weight
lattice in $\kt^*$. Then $\Lambda^*_+= \Lambda \cap \kt^*_+$ is the monoid of dominant weights.
Choose a set of fundamental weights $\Pi = \{\varpi_1,\cdots, \varpi_r \}$, which spans $\Lambda^*$ as
$\bZ$-basis. Let $V_{\varpi_i}$ be the irreducible representation of $K$ with the
highest weight $\varpi_i$ and $v_i$ be a fixed highest weight vector of $V_{\varpi_i}$.
We will show that $\uip$ can be embedded into $E=\oplus_{\varpi \in \Pi}V_{\varpi}$.

Before describing the embedding, we recall some facts about the Hamiltonian
action on the vector space $E$. Take a Hermitian metric $(-,-)_E$ on $E$ such
that one can decompose $E$ into direct sum of unitary $K$-subrepresentation
$V_{\varpi}$ and $\norm{v_i} = 1$. The symplectic form and moment map of $E$ are
given by
\begin{equation}
  \label{eq:vectsymp}
  \omega_E(v,w) = - \Im{{(v,w)}_E}\quad\text{and}\quad \langle{\mu_E(v)},X\rangle=\frac{1}{2}\omega_E(X.v,v),
\end{equation}
respectively, where $v, w\in E$ and $X\in \kk$. For the later usage, one also
defines a $T$-action on $E$ by requiring $T$ acting on $V_{\varpi}$ with the weight
$-\varpi$. Clearly, this $T$-action commutes with the $K$-action.

\subsection{}
\label{sub:liesymb}
Now let $\{\alpha_1,\cdots,\alpha_r\}\subseteq \kt^*$ be simple roots of $\kk$, and
$\{\alpha_1^{\vee},\cdots,\alpha_r^{\vee}\}\subseteq \kt$ be the corresponding coroots, i.e.\ $\alpha_i^{\vee} =
2\alpha_i^*/(\alpha_i,\alpha_i) \in \kt$, where $(-,-)$ is a Weyl group invariant inner product
on $\kt^*$ and $\alpha_i^*$ is the dual element of $\alpha_i$ with respect to the inner
product. It is well known that $\langle \varpi_i^{},\alpha_j^{\vee}\rangle = \delta_{ij}$ for $i,j\in
\{1,\cdots,r\}$. Since $\lambda(\alpha_i^{\vee})\ge 0$ holds for any $\lambda\in \kt^*_+$, we can define a
continuous map from $K \times \kt^*_+$ to $E$ as follows,
\begin{equation}
  \label{eq:deff}
  \cF(k,\lambda) = \frac{1}{\sqrt{\pi}}\sum^r_{p=1}\sqrt{\lambda(\alpha_p^{\vee})}k.v_p.
\end{equation}
By~\cite[Lemma~6.2]{Guillemin:2002sx}, $\cF$ can descend to a map on $\uip$ and
we will denote the descended map by the same symbol $\cF$. In fact, $\cF:\uip \rightarrow
E$ gives the claimed embedding. Let $G = K^{\bC}$ be the complexification of $K$
and $N$ be a maximal unipotent group. The image of $\cF$ in $E$ is
\begin{equation*}
  G_N \coloneq \overline{G(\sum_{\varpi\in \Pi}v_{\varpi})},
\end{equation*}
with respect to the Zariski or standard topology. As an affine subvariety of
$E$, $G_N$ has an induced Hamiltonian $K \times T$-space structure, which is
isomorphic to $\uip$ via $\cF$ by~\cite[Proposition~6,8]{Guillemin:2002sx}.
Therefore, by identified with $G_N$, $\uip$ is given a complex structure, which
is compactible with its symplectic structure. Using such an description of
$\uip$, one can find the principal cross-section is the unique open and dense
$G$ orbit $G(\sum_{\varpi\in \Pi}v_{\varpi})\simeq G/N$. Moreover, each stratum
in (\ref{eq:tkdc}) corresponds to a $G$ orbit, i.e.
\begin{equation*}
  \cF({K}/{[K_{\sigma},K_{\sigma}]} \times \sigma) = G/[P_{\sigma},P_{\sigma}],
\end{equation*}
where $P_{\sigma}$ is the parabolic subgroup associated to face $\sigma$.

\subsection{}
\label{sub:cpx-on-uip}
When $\uip$ is equipped with a complex structure, one can use (\ref{eq:up}) to
define the complex structure on $M_{\imp}$ as a symplectic quotient of a \kah\
manifold in the usual way which we review in the next section. A small
inaccuracy here is that since $M \times \uip$ is a \kah\ space, not smooth
necessarily, the symplectic reduction needs some extension, that is, one should
use the reduction of complex spaces developed
in~\cite{Heinzner:1994bn,Heinzner:1999aa}. However, for our purpose, we actually
don't need such a general theory. The reason is that the space we are interested
can always be obtained via the symplectic reduction construction involving only
a suitable smooth stratum of $\uip$.

As we have recalled, the $T$-action on $M_{\imp}$ is Hamiltonian. In the \kah\
case, one can further assert that the $T$-action is also holomorphic by the
result in~\cite{Greb_2010ka}. We give a proof for the simpler manifold case in
the next section.

\section{The \nth{1} proof of Theorem~\ref{thm:intro1}: the symplectic implosion approach}
\label{sec:1st-app}

\subsection{}
In this section, we use the \kah\ geometry properties of the symplectic
implosion to prove Theorem~\ref{thm:intro1}. To this aim, we first give a proof
for Theorem~\ref{thm:introeq} by using Proposition~\ref{prop:red-in-stage},
which is a ``reduction in stages'' result for \kah\ manifolds. Then, we show
Theorem~\ref{thm:intro1} by combining Theorem~\ref{thm:introeq} and the property
of the abelian reduction.

\subsection{\kah\ structure on symplectic quotients.}
\label{sub:kahstr}
For readers' convenience, we briefly recall the definition of the \kah\
structure on a symplectic quotient. Let $(X,\omega_X,\cst_X,\mu_X)$ be a Hamiltonian
\kah\ $K$-manifold. Suppose that $0$ is a regular value of the moment map
$\mu_X$.\footnote{In general, any regular value lying in $\kz^*$ works, where $\kz$
  is the center of $\kk$.} At any point $x\in \mu^{-1}(0)$, one has the following
orthogonal decomposition of the tangent space.
\begin{equation}
  \label{eq:tangdecomp}
  \rT_x \mu^{-1}(0) = Q_x \oplus \kk\cdot x \quad \text{and} \quad \rT_x X = Q_x \oplus \kk \cdot x \oplus \cst_X {(\kk \cdot x)},
\end{equation}
where $\kk \cdot x$ is tangent subspace induced by the infinitesimal $K$-action at
$x$. By the definition of a \kah\ metric, $Q_x$ is a $\cst$-invariant subspace,
which implies that $Q_x$ is also a symplectic subspace of $\rT_x X$. Hence,
\begin{equation*}
  \diff \pi_x:Q_x \rightarrow \rT_{\pi(x)} X_0
\end{equation*}
is an isomorphism, where $\pi$ is the quotient
map from $\mu^{-1}(0)$ to the symplectic quotient $X_0= \mu^{-1}(0)/K$. Via $\diff
\pi(x)$, the (Riemannian) metric and the almost complex structure on $\rT_{\pi(x)}
X_0$ are induced from those of $Q_x$. One can show that,
\cite[Lemma~7.2.7]{Futaki:1988aa} or \cite[Lemma~2.2]{La-Nave:2016en}, the
almost complex structure on $X_0$ obtained in this way is integrable and
compatible with the reduced symplectic form $\omega_0$ on $X_0$.

As an additional remark, to define a complex structure on a symplectic quotient,
besides the method described as above, there is another possible way, i.e.\ by
selecting a semi-stable point set, $X^{\rss}$, of $X$ and defining $X_0$ to be
the holomorphic quotient of $X^{\rss}$ by $K^{\bC}$,
c.f.~\cite[\S~7]{Kirwan:1984st}. Readers may find that a short review of this method
in \S~\ref{sec:ss} is helpful.

Now we show that the $T$-action on a symplectic implosion is holomorphic. Since
the complex structure on a symplectic implosion is defined using symplectic
reduction, this result follows from the following general property about group
actions on a symplectic quotient.

Before the proof, we remark that the following two propositions are well known
to experts. In fact, as the referee points to us, similar results are even
proved for complex spaces in~\cite{Greb_2010ka}. However, since our manifold
settings are much simpler than~\cite{Greb_2010ka}, we would like to provide
shorter differential geometry proofs for these results to make this paper more
self-contained.

\begin{proposition}
  \label{prop:induced-action}
  Assume that there is another compact, Hamiltonian and holomorphic group action
  $L$ on $X$. Moreover, the $L$-action commutes with the $K$-action and
  preserves the moment map $\mu$ of the $K$-action. Then $L$ induces an action on
  the symplectic quotient $X_0$, which preserves the \kah\ structure on $X_0$.
\end{proposition}

\begin{proof}
  It is a classical result that the induced $L$-action preserves the reduced
  symplectic structure, c.f.~\cite[Theorem~2]{Marsden:1974aa}. We are going to
  show the induced $L$-action also preserves the complex structure on $X_0$. Let
  $W\in \mathrm{Lie}(L)$ and $W^X$ be the induced vector field of $W$ on
  $X$.\footnote{The restriction of $W^X$ on $\mu^{-1}(0)$ is denoted by the same
    symbol.} First, we show that the Lie derivative of $W^X$ preserves sections
  lying in the subbundle $Q\subseteq \rT \mu^{-1}(0)$ over $\mu^{-1}(0)$. For any $\bfw\in
  \Gamma(Q)$ being a smooth section of $Q$, if we view $\bfw$ as a vector field of
  the submanifold $\mu^{-1}(0)$,
  \begin{equation*}
    \sL_{W^X}{\bfw} = [W^X,\bfw] \in \Gamma(\rT \mu^{-1}(0))
  \end{equation*}
  is well-defined. For any $Y\in \kk$, since $L$ preserves the restricted metric
  on $\mu^{-1}(0)$, one has
  \begin{equation*}
    (\sL_{W^X}{\bfw}, Y^X) = W^X(\bfw,Y^X) - (\bfw, [W^X,Y^X]) = 0,
  \end{equation*}
  which implies that $\sL_{W^X}{\bfw}\in \Gamma(Q)$.

  Now choose any vector field $\bfv$ on $X_0$, denoting the lift-up of $\bfv$ in
  $\Gamma(Q)$ by $\bfv^{\sharp}$. Let $\cst_0$ be the reduced complex structure on $M_0$.
  By the definition of $\cst_0$ and the induced $L$-action, one has
  \begin{equation}
    \label{eq:red-cpx-def}
    \cst_0{\bfv} = \diff \pi(\cst{\bfv^{\sharp}}),\quad W^{X_0} = \diff \pi (W^{X}).
  \end{equation}
  Using (\ref{eq:red-cpx-def}), we can calculate as follows,
  \begin{equation}
    \label{eq:red-cpx-str}
    \begin{aligned}
      (\sL_{W^{X_0}}{\cst_0})(\bfv) & = [{W^{X_0}},\cst_0(\bfv)] -
      \cst_0{[W^{X_0},\bfv]}\\
      & = [\diff \pi(W^X), \diff \pi(\cst{\bfv^{\sharp}})] - \cst_0 [\diff \pi(W^X),\diff
      \pi(\bfv^{\sharp})] \\
      & = \diff \pi ([W^X, \cst{\bfv^{\sharp}}]) - \cst_0 \diff \pi ([W^X, \bfv^{\sharp}])\\
      & = \diff \pi ([W^X, \cst{\bfv^{\sharp}}]) - \diff \pi \cst ([W^X, \bfv^{\sharp}])\\
      & = \diff \pi ((\sL_{W^X}{\cst})(\bfv^{\sharp})) = 0.
    \end{aligned}
  \end{equation}
  Note that in the \nth{4} equality, we have used $[W^X, \bfv^{\sharp}] =
  \sL_{W^X}{\bfv^{\sharp}}\in \Gamma(Q)$ and $[W^X, \bfv^{\sharp}]$ is a $K$-invariant vector
  field. Since $\bfv$ is arbitrary, one concludes that
  $\sL_{W^{X_0}}{\cst_0}=0$.
\end{proof}

With Proposition~\ref{prop:induced-action}, we can state the following ``reduction
in stages'' result for \kah\ manifolds, which is a preparation for the proof of
Theorem~\ref{thm:introeq}.
\begin{proposition}
  \label{prop:red-in-stage}
  Let $(X,\omega_X,\cst_X,\mu_1\times \mu_2)$ be a Hamiltonian \kah\ $K_1 \times T_2$-manifold ,
  where $K_1$ is a compact semi-simple group and $T_2$ is a compact torus.
  Suppose that $(0,\lambda)$ is a regular value of $\mu_1\times \mu_2$. Then as \kah\
  manifolds, $X\dsl_{(0,\lambda)} (K_1\times T_2)$ is naturally isomorphic to $(X\dsl_0
  K_1) \dsl_{\lambda} T_2$.
\end{proposition}

The two symplectic quotients in the proposition are canonically diffeomorphic to
each other and we will identify the two quotients in the following. By using
reduction in stages for symplectic manifolds (or spaces),
c.f.~\cite[Theorem~2]{Marsden:1974aa} and~\cite[\S~4]{Sjamaar:1991aa}, the
reduced symplectic forms given by the two methods in the proposition are the
same. Therefore, to show the two reduction procedures lead to the same complex
structure on the symplectic quotient, one only needs to show that the two
reduced metrics on the symplectic quotient coincide. For this purpose, the
following elementary fact about projections of linear spaces is useful.

\begin{lemma}
  \label{lm:proj}
  Let $E$ be a complex vector space with a Hermitian inner product. $V$ and $W$
  are two subspaces of $E$. For the following three orthogonal projections,
  $\prj:E \rightarrow (V + W)^{\perp}$, $\prj_1:E \rightarrow W^{\perp}$, and $\prj_2:W^{\perp} \rightarrow
  \big(\prj_1(V)\big)^{\perp}\cap W^{\perp}$, one has $\prj = \prj_2\prj_1$.
\end{lemma}

\begin{proof}
  If $u \in V + W$, one has $u = v + w$ with $v \in V$, $w \in W$. Then by
  definitions, $\prj_2\prj_1(u) = \prj_2 \prj_1(v) = 0$, and $\prj(u) = 0$. If
  $u \in (V + W)^{\perp}$, one has $u \in W^{\perp} \cap V^{\perp}$. Therefore, $\prj_2\prj_1(u) =
  \prj_2(u)$. On the other hand, if $u\in (V + W)^{\perp}$, since any $v' \in \prj_1(V)$
  can be written as $v' = v - v''$ with $v \in V$, $v''\in W$, one has $u \perp
  v'$, i.e. $u\in \big(\prj_1(V)\big)^{\perp}$. Hence, $\prj_2(u) = u = \prj(u)$.
\end{proof}

\begin{proof}[Proof of Proposition~\ref{prop:red-in-stage}]
  Firstly, we assume $X$ to be a manifold. Let $x\in X$ lying in the level set
  $(\mu_1\times \mu_2)^{-1}(0,\lambda)$. Denote $W$ to be the subspace of $\rT_{x}X$ generated
  by the infinitesimal $K_1$-action and $V\subseteq \rT_xX$ to be the subspace generated
  by the infinitesimal $T_2$-action. Choose
  \begin{equation*}
    \bfv \in \rT_x(\mu_1\times \mu_2)^{-1}(0,\lambda)\subseteq \rT_xX
  \end{equation*}
  transversal to $W+ V$. Introduce the following two quotient maps,
  \begin{align*}
    \pi:&(\mu_1\times \mu_2)^{-1}(0,\lambda) \rightarrow (\mu_1\times \mu_2)^{-1}(0,\lambda)/(K_1\times T_2) = X
    \dsl_{(0,\lambda)}(K_1\times T_2),\\
    \bar{\pi}:&(\mu_1\times \mu_2)^{-1}(0,\lambda) \rightarrow \big((\mu_1\times \mu_2)^{-1}(0,\lambda)/K_1\big)/T_2 = (X \dsl_{0}K_1)
    \dsl_{\lambda} T_2.
  \end{align*}
  Since $(\mu_1\times \mu_2)^{-1}(0,\lambda)/(K_1\times T_2)$ is naturally isomorphic to $\big((\mu_1\times
  \mu_2)^{-1}(0,\lambda)/K_1\big)/T_2$, $\pi$ and $\bar{\pi}$ are essentially the same map.
  We decide to use different notations to remind us that they originate from
  different reduction procedures. We are going to show that the norm on $\diff
  \pi(\bfv)$, coming from reduction of $K_1\times T_2$, is equal to the norm on $\diff
  \bar{\pi}(\bfv)$, coming from reduction first by $K_1$ and then by $T_2$. As a
  result, the two metrics on the common quotient coincide. Suppose $\prj$ to be
  the orthogonal projection from $\rT_x (\mu_1\times \mu_2)^{-1}(0,\lambda)$ onto $(V +
  W)^{\perp}$, $\prj_1$ to be the orthogonal projection from $\rT_x(\mu_1\times
  \mu_2)^{-1}(0,\lambda)$ onto $W^{\perp}$ and $\prj_2$ to be the orthogonal projection from
  $W^{\perp}$ onto $\big(\prj_1(V)\big)^{\perp} \cap W^{\perp}$. By the definition of the
  metric on a symplectic quotient, the norm of $\diff \pi(\bfv)$ is equal to
  $\norm{\prj(\bfv)}$ and the norm of $\diff \bar{\pi}(\bfv)$ is equal to
  $\norm{\prj_2\prj_1(\bfv)}$. Therefore, the equality of the metrics is a
  result of Lemma~\ref{lm:proj}.
 
\end{proof}

\begin{remark} 
  \label{rk:red-in-stage}
  We make some comments about implications of Proposition~\ref{prop:red-in-stage}.
  \begin{enumerate}[label=\Alph*.]
  \item\label{rk:red-a} In Proposition~\ref{prop:red-in-stage}, the $K_1$-action and the
    $T_2$-action play similar roles, which means that it makes no difference that
    which action comes first when we perform reduction in stages. As a corollary,
    one can see that
    \begin{equation*}
      (X\dsl_{0} K_1)\dsl_{\lambda} T_2 \simeq (X\dsl_{\lambda} T_2)\dsl_{0} K_1.
    \end{equation*}
\item The same proof of Proposition~\ref{prop:red-in-stage} also works for a
    more general case: $K_1\times T_2$ can be replaced by $K_1\times K_2$, where $K_1,K_2$
    are arbitrary compact groups; $(0,\lambda)$ can be replaced by $(\lambda_1,\lambda_2)$, where
    $\lambda_i$ lies in the dual of the center of $\kk_i$, $i=1,2$.
\end{enumerate}
\end{remark}

\subsection{}
Let $(M,\omega,\cst,\mu)$ be a Hamiltonian \kah\ $K$-manifold.
As an application of Proposition~\ref{prop:red-in-stage}, we first give a proof of
Theorem~\ref{thm:introeq}.

\begin{proof}[Proof of Theorem~\ref{thm:introeq}]
  We first prove the case that $K$ is a semi-simple simply connected group.
  We begin with an analysis of the metric of $\uip$. Here, we use notational
  conventions as in \S~\ref{sub:liesymb}. Let $(R,R_+)$ be a root system for
  $\kg = \kk^{\bC}$ containing simple roots $S=\{\alpha_1,\cdots,\alpha_r\}$.
  One has the following decomposition of $\kk$,
  \begin{equation}
    \label{eq:kkdecomp}
    \kk = \kt \oplus \sum_{\alpha \in R_+}\kr_{\alpha},\quad\kt =  \sum_{i=1}^{r}\bR \cdot \alpha_i^{\vee},
  \end{equation}
  where $\kr_{\alpha} =\big((\kg)_{\alpha}\oplus (\kg)_{-\alpha}\big)\cap \kk$
  is a real two-dimensional vector space. Let $R(\sigma)\subseteq R$ be a subset
  of roots such that $\langle\xi,\alpha^{\vee}\rangle=0$ for any $\xi\in \sigma,
  \alpha\in R(\sigma)$ and let $S(\sigma) \coloneq S \cap R(\sigma)$. Recall
  that $\kk_{\sigma}\subseteq \kk$ is the Lie algebra of the isotropic group of
  any point lying in $\sigma$. The following decomposition also holds.
  \begin{equation}
    \label{eq:ksigmadecomp}
    \begin{aligned}
      [\kk_{\sigma},\kk_{\sigma}] &=& \sum_{\alpha_i\in S(\sigma)}\bR \cdot \alpha_i^{\vee} \oplus \sum_{\alpha \in R(\sigma)\cap R_+}\kr_{\alpha}, \\
      \kk/[\kk_{\sigma},\kk_{\sigma}] &\simeq& \sum_{\alpha_i\notin
        S(\sigma)}\bR \cdot \alpha_i^{\vee} \oplus \sum_{\alpha \in R_+\setminus
        R(\sigma)}\kr_{\alpha}.
    \end{aligned}
  \end{equation}
  One notices that for the smooth stratum $K/[K_{\sigma},K_{\sigma}]\times
  \sigma$ of $\uip$, the tangent space at $([e],\lambda)$ is
  $\kk/[\kk_{\sigma},\kk_{\sigma}]\oplus \sigma$. By (\ref{eq:kmm}), for any
  $([k],\lambda)\in \uip$, one has $\Phi_{\imp,\cR}([k],\lambda) = \lambda$,
  which means that
  \begin{equation*}
    \rT_{([e],\lambda)}\Phi_{\imp,\cR}^{-1}(\lambda) = \kk/[\kk_{\sigma},\kk_{\sigma}].
  \end{equation*}

  Recall that in \S~\ref{sub:liesymb}, we define the complex structure, or the
  \kah\ structure equivalently, on $\uip$ using an embedding $\cF:\uip
  \rightarrow E$. By the definition of the symplectic form on $E$,
  (\ref{eq:vectsymp}), for any $H\in \sum_{\alpha_i\notin S(\sigma)}\bR \cdot
  \alpha_i^{\vee}$ or $X\in \kr_{\alpha}$, $\alpha\in R_+\setminus R(\sigma)$,
  one has
  \begin{equation}
    \label{eq:diff-of-f}
    \diff \cF_{([e],\lambda)}(V) = \frac{1}{\sqrt{\pi}}\sum^r_{p=1}\sqrt{\lambda(\alpha_p^{\vee})}V.v_p,\quad\text{for }V=H\text{ or }X.
\end{equation}
  Since $\{v_p\}$ are the highest weight vectors, the above equality implies
  that $(H.v_p,X.v_q)=0$ for any $p,q\in \{1,\cdots,r\}$. Therefore,
  $\sum_{\alpha_i\notin S(\sigma)}\bR \cdot \alpha_i^{\vee}$ and $\sum_{\alpha
    \in R_+\setminus R(\sigma)}\kr_{\alpha}$ as subspaces at
  $\rT_{([e],\lambda)} \Phi_{\imp,\cR}^{-1}(\lambda)$ are always orthogonal to
  each other under the pullback metric. Clearly, the former subspace is the
  subspace generated by the right $T$-action at this point. To calculate the
  metric on the symplectic quotient of the $T$-action, we only need to calculate
  the metric on $\sum_{\alpha \in R_+\setminus R(\sigma)}\kr_{\alpha}\subseteq
  \kk/[\kk_{\sigma},\kk_{\sigma}]$.
  
  Choose a vector $X_{\alpha}$ in $\kg_{\alpha}$ for $\alpha \in R_+\setminus
  R(\sigma)$. Let $\theta$ be the Cartan involution corresponding to
  $\kk\subseteq \kg$. One has $X_{-\alpha}\coloneq \theta(X_{\alpha})\in
  \kg_{-\alpha}$. Using a suitable normalization of $X_{\alpha}$, we can assume
  that $[X_{-\alpha},X_{\alpha}]=-\opi \alpha^{\vee}$. By the definition of the
  Cartan involution, $U_{\alpha} \coloneq X_{\alpha} + X_{-\alpha}$ and
  $V_{\alpha} \coloneq \opi X_{\alpha} - \opi X_{-\alpha} $ are vectors in
  $\kr_{\alpha}$. For any $\alpha,\alpha'\in R_+\setminus R(\sigma)$, using
  (\ref{eq:vectsymp}) and (\ref{eq:diff-of-f}), we can calculate the inner
  product of $U_{\alpha}$ and $V_{\alpha'}$ as vectors in
  $\rT_{([e],\lambda)}\Phi_{\imp,\cR}^{-1}(\lambda)$ under the pullback metric
  of $E$ in the following way, \begingroup \allowdisplaybreaks
  \begin{align*}
(U_{\alpha},V_{\alpha'})_{([e],\lambda)}
    & = \frac{-1}{\pi} \sum^r_{p=1}{\lambda(\alpha_p^{\vee})}\Im{(U_{\alpha}.v_p,\opi V_{\alpha'}.v_p)}_E \\
    &= \frac{-1}{\pi} \sum^r_{p=1}{\lambda(\alpha_p^{\vee})} \Im{(X_{-\alpha}.v_p, X_{-\alpha'}.v_p)}_E \\
    & = \frac{-1}{\pi} \sum^r_{p=1}{\lambda(\alpha_p^{\vee})} \Im{(-X_{\alpha'}.X_{-\alpha}.v_p, v_p)}_E \\
    &=  \frac{-1}{\pi} \sum^r_{p=1}{\lambda(\alpha_p^{\vee})} \Im{([X_{\alpha'},X_{-\alpha}].v_p, v_p)}_E.
  \end{align*}
  \endgroup In the above equalities, the \nth{1} one uses the relation between
  the real inner product and the Hermitian inner product on $E$. The \nth{2} and
  \nth{4} equalities use the fact that $v_p$ is the highest weight vector. As
  for the \nth{3} equality, recall that the adjoint operator of $X_{-\alpha}$ is
  $-X_{\alpha}$.
  
  Now if $\alpha \neq \alpha'$, then $[X_{\alpha'},X_{-\alpha}]\notin \kt$,
  consequently $([X_{\alpha'},X_{-\alpha}].v_p,v_p)_E=0$. Otherwise, for $\alpha
  = \alpha'$, one has
  \begin{equation*}
    \Im{([X_{\alpha'},X_{-\alpha}].v_p, v_p)}_E = \Im{(\opi \alpha^{\vee}.v_p, v_p)}_E = -\Im{(2\pi \varpi_p(\alpha^{\vee}))} = 0.
  \end{equation*}
  All in all, the inner product of $U_{\alpha}$ and $V_{\alpha'}$ always
  vanishes. In the same way, one can show that
  \begin{equation}
    \label{eq:inner-prod-uv}
    (U_{\alpha},U_{\alpha'})_{([e],\lambda)} =
    \begin{cases}
      2 \sum^r_{p=1}{\lambda(\alpha_p^{\vee})} \varpi_p(\alpha^{\vee}) = 2 \langle\lambda,\alpha^{\vee}\rangle &,\;  \alpha=\alpha'; \\
      0 &,\; \alpha \neq \alpha',
    \end{cases}
  \end{equation}
  where we have used $\langle\varpi^{}_p,\alpha^{\vee}_q\rangle = \delta_{pq}$.
  For the inner product between $V_{\alpha}$ and $V_{\alpha'}$, the result is
  similar.

  Recall that the Kirillov-Kostant-Souriau symplectic form on the coadjoint
  orbit $\cO_{\lambda} \coloneq K\cdot\lambda \subseteq \kk^*$ is defined as
  follows,
  \begin{equation}
    \label{eq:kks}
    \omega_{\cO_{\lambda}}(X\cdot[\lambda],Y\cdot[\lambda])  = \langle\lambda,[X,Y]\rangle,
  \end{equation}
  where $X,Y\in \kk$. As we have said before, the symplectic reduction of $\uip$
  at $\lambda$ for the $T$-action is $\cO_{\lambda}^*$, i.e. $\cO_{\lambda}$
  with the symplectic form
  $-\omega_{\cO_{\lambda}}$.\footnote{$\cO_{\lambda}^*$, as a symplectic
    manifold, is isomorphic to the coadjoint orbit $\cO_{-\mathrm{w}_0
      \lambda}$, where $\mathrm{w}_0$ is the longest element of the Weyl group.}
  Now if we define the almost complex structure on $\rT_{\lambda}\cO_{\lambda}^*
  \simeq \sum_{\alpha \in R_+\setminus R(\sigma)}\kr_{\alpha}$ as follows, for
  $\alpha\in R_+\setminus R(\sigma)$,
  \begin{equation}
    \label{eq:cst-on-orbit}
    \cst {(U_{\alpha}\cdot[\lambda])} = -V_{\alpha}\cdot[\lambda],\quad \cst{(V_{a}\cdot[\lambda])} = U_{\alpha}\cdot[\lambda],
  \end{equation}
  it is well known that such an almost complex structure is integrable on the
  coadjoint orbit. By (\ref{eq:kks}) and (\ref{eq:cst-on-orbit}), the \kah\
  metric on $\cO_{\lambda}^*$ is,
  \begin{multline}
    \label{eq:metric-on-orbit}
    (U_{\alpha}\cdot[\lambda],U_{\alpha'}\cdot[\lambda])_{\cO_{\lambda}^*} =
    -\omega_{\cO_{\lambda}}(U_{\alpha}\cdot[\lambda],\cst{(U_{\alpha'}\cdot[\lambda])}) \\ = \langle\lambda,[U_{\alpha},V_{\alpha'}]\rangle
    = \begin{cases}
      2 \langle\lambda,\alpha^{\vee}\rangle &,\;  \alpha=\alpha'; \\
      0 &,\; \alpha \neq \alpha'.
      \end{cases}
  \end{multline}
  Terms involving $V_{\alpha}$ can be calculated in the same way. Comparing
  (\ref{eq:inner-prod-uv}) and (\ref{eq:metric-on-orbit}), one can conclude
  \begin{equation*}
    \uip \dsl_{\lambda} T =( K/[K_{\sigma},K_{\sigma}]\times \sigma) \dsl_{\lambda} T \simeq \cO_{\lambda}^*
  \end{equation*}
  as the \kah\ manifold. Now, by the
  definition of the \kah\ reduction, Proposition~\ref{prop:red-in-stage} and the
  item~\ref{rk:red-a} in Remark~\ref{rk:red-in-stage}, we can show the
  conclusion of Theorem~\ref{thm:introeq} using the following argument,
  \begin{multline}
    \label{eq:sqeq}
    M \dsl_{\lambda} K = (M \times \cO_{\lambda}^*)\dsl K \simeq \big((M \times
    K/[K_{\sigma},K_{\sigma}]\times
    \sigma) \dsl_{\lambda} T\big)\dsl K \\
    \simeq \big((M \times K/[K_{\sigma},K_{\sigma}]\times \sigma)\dsl
    K\big)\dsl_{\lambda} T = M_{\imp}\dsl_{\lambda} T,
  \end{multline}
  where in the last equality we use the fact that
  \begin{equation*}
    \Phi_{\imp,\cR}^{-1}(\lambda)\subseteq K/[K_{\sigma},K_{\sigma}]\times \sigma.
  \end{equation*}

  For a general compact group $K$, let $K = K_{\rss}Z$, where $K_{\rss}= [K,K]$
  is the semi-simple part of $K$ and $Z$ is the center of $K$. Choose a Cartan
  subgroup $T_{\rss}$ of $K_{\rss}$, then $T=T_{\rss}Z$ is a Cartan subgroup of
  $K$. Take $\lambda = \lambda_1 + \lambda_2 \in \kt^*_{\rss} \oplus \kz^* = \kt^*$, where $\lambda_1$ lies in
  the face $\sigma$ of $\kt^*_+$. Then $\lambda$ lies in the face $\bar{\sigma} \coloneq \sigma \oplus \kz^*$ of
  the positive Weyl chamber of $K$. Recall that in \S~\ref{sub:group-ass}, for
  a general $K$, the complex structure on $M_{\imp,K}$ is defined by
  $M_{\imp,K_{\rss}}$. Then, by Proposition~\ref{prop:red-in-stage} and
  (\ref{eq:sqeq}),
  \begin{align*}
    M \dsl_{\lambda}K & \simeq (M \times \cO_{\lambda_1}^*)\dsl_{(0,\lambda_2)}(K_{\rss} \times Z) \simeq \big((M \times \cO_{\lambda_1}^*)\dsl K_{\rss}\big)\dsl_{\lambda_2} Z \\
                & \simeq \big((M \times K_{\rss}/[K_{\rss,\sigma},K_{\rss,\sigma}] \times \sigma)\dsl K_{\rss}\big)\dsl_{\lambda_1} T_{\rss}\dsl_{\lambda_2} Z \\
                & \simeq \big((M \times (\rT^* K_{\rss})_{\imp})\dsl K_{\rss}\big) \dsl_{\lambda} T \simeq M_{\imp,K} \dsl_{\lambda} T,
\end{align*}
  where $\cO_{\lambda_1}$ is the coadjoint orbit of $\lambda_1$ in ${\kk_{\rss}}$.
\end{proof}

As an application of Theorem~\ref{thm:introeq}, we can show the local invariance
of complex structures on a symplectic quotient, i.e.~Theorem~\ref{thm:intro1}.
Firstly, we recall the corresponding result for a torus action. For simplicity,
we only state the result for an $S^1$-action. The higher dimensional case is
similar.

\begin{lemma}
  \label{lm:ab-eq}
  Let $K=S^1$. Instead of assuming the compactness of $M$, here we only require
  the moment map $\mu:M\rightarrow\kk^*\simeq \bR^1$ of the $K$-action is
  proper. If $c$ is a regular value of $\mu$ and $S^1$ acts on $\mu^{-1}(c)$
  freely, then for any $c'\in \bR$ close to $c$ enough, $M_{c'}$ is
  biholomorphic to $M_{c}$.
\end{lemma}
For readers' convenience, we reproduce the proof
of~\cite[Lemma~7.4.1]{Futaki:1988aa}. One can also
see~\cite[Proposition~2.4]{La-Nave:2016en}.

\begin{proof}
  Choosing a non-zero vector $X$ in $\kk$, since $c$ is a regular value of $\mu$
  and $c,c'$ are close enough, the flow generated by $\cst{X^M}$ induces a
  diffeomorphism between $M_{c'}$ and $M_{c}$. We will show that this
  diffeomorphism is also holomorphic.

  Since the almost complex structure on $M$ is integrable, the corresponding
  Nijenhuis tensor vanishes. Therefore, for any vector field $\bfv$,
  \begin{equation}
    \label{eq:nij}
    [\cst{X^M},\cst{\bfv}] - \cst{[\cst{X^M},\bfv]} = \cst{[X^M,\cst{\bfv}]} + [X^M,\bfv].
  \end{equation}
  Since $X^M$ induces a holomorphic isometry, one has
  \begin{equation}
    \label{eq:hol-iso}
    0 = (\sL_{X^M}\cst)({\bfv}) = \sL_{X^M}(\cst{\bfv}) - \cst{(\sL_{X^M} \bfv)}
    = [X^M,\cst{\bfv}] - \cst{[X^M,\bfv]}.
  \end{equation}
  Combine (\ref{eq:nij}) and (\ref{eq:hol-iso}),
  \begin{equation}
    \label{eq:jx-hol}
    (\sL_{\cst{X^M}}\cst)({\bfv}) = \sL_{\cst{X^M}}(\cst{\bfv}) - \cst{(\sL_{\cst{X^M}} \bfv)}
    = [\cst{X^M},\cst{\bfv}] - \cst{[\cst{X^M},\bfv]} = 0.
  \end{equation}
  Since $\bfv$ is an arbitrary vector field, (\ref{eq:jx-hol}) implies that the
  flow generated by $\cst{X^M}$ preserves the complex structure of $M$. Note
  that this flow also preserves the decomposition of (\ref{eq:tangdecomp}).
  Therefore, the complex structure of the subspace $Q$ is invariant under the
  action of the flow, which implies that the diffeomorphism between $M_{c}$ and
  $M_{c'}$ is holomorphic.

\end{proof}

\subsection{}
Due to Theorem~\ref{thm:introeq}, Theorem~\ref{thm:intro1} is an easy corollary
of Lemma~\ref{lm:ab-eq}. More precisely, we use a minor generalization of
Lemma~\ref{lm:ab-eq}: instead of assuming the properness of $\mu$, we only need
that $\mu$ is proper in a neighborhood of $\mu^{-1}(c)$, under which condition
the proof of the lemma still works.
\begin{proof}[Proof of Theorem \ref{thm:intro1}]

  Using Proposition~\ref{prop:red-in-stage} and Lemma~\ref{lm:ab-eq}, we can
  see the only case needed to show is that $K$ is semi-simple and
  simply-connected. We choose a small neighborhood $U$ of $\lambda$ in $\sigma$ and denote
  $S\coloneq (M \times K/[K_{\sigma},K_{\sigma}]\times U)\dsl K$. Since $\lambda$ is regular, $S$ is a smooth
  manifold. For $\lambda,\lambda'\in U$, by the proof of Theorem~\ref{thm:introeq}, we have
  $M_{\lambda} = S \dsl_{\lambda} T$, $M_{\lambda'} = S \dsl_{\lambda'} T$. Now, by Lemma~\ref{lm:ab-eq}
  and the remark in the previous paragraph, $S \dsl_{\lambda} T$ is biholomorphic to
  $S \dsl_{\lambda'} T$, which implies the result of Theorem \ref{thm:intro1}.

\end{proof}

\begin{remark*}
  If we only want to show Theorem~\ref{thm:intro1}, from the above proof, the
  introduction of the symplectic implosion is not so necessary. In fact, we only
  use a suitable smooth stratum of the symplectic implosion. This mainly because
  our problem is local in essence. However, if we remove all the material about
  the symplectic implosion, the whole proof would be in an ad hoc flavor more or
  less. In our humble opinion, the symplectic implosion gives the right
  framework to understand the problem.

\end{remark*}

\section{A GIT explanation of Theorem~\ref{thm:intro1}}
\label{sec:ss}

\subsection{}
It is a seminal result of~\cite{Kirwan:1984st,Ness:1984rr} that a symplectic
quotient is naturally isomorphic to a GIT quotient. Inspired by this fact, we
discuss a result in this section, Proposition~\ref{prop:ss}, which can be seen
as a reinterpretation of Theorem~\ref{thm:intro1} in the GIT language.

\subsection{Stability conditions.}
\label{sub:ss-point}
Firstly, let us recall the definition of a key concept: semi-stable (or stable) points. As in \S~\ref{sub:liesymb}, let $G = K^{\bC}$ be the
complexifiction of the compact group $K$. For a compact Hamiltonian \kah\
$K$-manifold $(M,\omega,\mu,\cst)$, there is a holomorphic $G$-action on $M$ induced by
the $K$-action. Any point $m\in M$ is called a \emph{semi-stable} point, if
$\overline{G\cdot m} \cap \mu^{-1}(0)\neq \emptyset$. Furthermore, if $m$ also satisfies that $G\cdot m
\cap \mu^{-1}(0)\neq \emptyset$ and the isotropic group $G_m$ at $m$ is finite, $m$ is called a
\emph{stable} point. We denote the sets of semi-stable and stable points of $M$
by $M^{\rss}$ and $M^{\rs}$ respectively. Sometimes we also use the notation
$M^{\rss}(\mu)$ and $M^{\rs}(\mu)$ to emphasize the dependence of semi-stable and
stable point sets on the moment map.

\begin{remark}
  \label{rk:ss-def}
  A variety of equivalent properties has been used to define the stability
  condition. Especially for the analytic stability condition on a \kah\
  manifold, which we use here, the definition may vary from paper to paper in
  literature. Therefore, we would like to make some clarifications about the
  terminology used in this paper. In short, the definition of semi-stable points
  in \S~\ref{sub:ss-point} coincides with the definition of semi-stable points
  appeared in~\cite{Heinzner:1994bn,Heinzner:1999aa}. It turns out that the set
  $M^{\rs}$ has appeared in~\cite{Kirwan:1984st} under the name the minimal
  stratum\footnote{Note that in~\cite{Sjamaar:1995aa}, points in this set are
    called analytic semi-stable.}, which is defined using the gradient flow of
  $\norm{\mu}^2$. For a proof for the coincidence of $M^{\rs}$ and the minimal
  stratum, one can see~\cite[Theorem~7.4]{Kirwan:1984st}
  and~\cite[Proposition~2.4]{Sjamaar:1995aa}.

  We also remind readers that the terminology of semi-stable points also appears
  in~\cite{Kirwan:1984st}. But it seems that the author reserves this concept
  for algebraic manifolds exclusively therein and the definition conforms to the
  stardard GIT one,~\cite{Mumford:1994aa}.
\end{remark}

\subsection{}
\label{sub:kirss}
Since we always assume $0$ to be a regular value of the moment map and $K$ to
act on $\mu^{-1}(0)$ freely, the semi-stable or stable point set behaves
particularly well. In fact, by~\cite[Theorem~7.4 {\&} 8.10]{Kirwan:1984st},
$M^{\rss}$ and $M^{\rs}$ coincide in this case. In particular, $M^{\rss}$ is a
$G$ principal bundle and one can define the complex structure on the symplectic
quotient by the formula: $M^{\rss}/G=\mu^{-1}(0)/K=M_0$. As we have said, Kirwan
and Ness's theorem ensures that such a definition of the complex structure on a
symplectic quotient is the same as the definition given in \S~\ref{sub:kahstr}.

Recall that $M_{\lambda}$, $M_{\lambda'}$ can be viewed as two symplectic
quotients of a common manifold $M \times K/K_{\sigma}$ with respect to two
different symplectic structures. In view of the above equivalent definition of
the complex structure on a symplectic quotient, we can explain
Theorem~\ref{thm:intro1} as follows. The two different symplectic structures on
$M \times K/K_{\sigma}$ give the same semi-stable point set. More precisely, one
has the following result.

\begin{proposition}
  \label{prop:ss}
  Let $\lambda,\lambda'$ be two points lying in a face $\sigma\subseteq
  \kt^*_+$. Identifying $K/K_{\sigma}$ with $\cO^*_{\lambda}$ using the map
  $[k]\mapsto k\cdot\lambda$, one obtains a Hamiltonian \kah\ structure on $M
  \times K/K_{\sigma}$, whose moment map is denoted by $\Xi_{\lambda}$.
  Similarly, by substituting $\cO^*_{\lambda}$ with $\cO^*_{\lambda'}$, one can
  define another moment map $\Xi_{\lambda'}$ on $M\times K/K_{\sigma}$. If
  $\lambda,\lambda'$ are close enough, the following two semi-stable point sets
  are equal:
  \begin{equation*}
    (M\times K/K_{\sigma})^{\rss}(\Xi_{\lambda}) = (M\times K/K_{\sigma})^{\rss}(\Xi_{\lambda'}).
  \end{equation*}
\end{proposition}

It is possible to prove Proposition~\ref{prop:ss} using techniques from GIT
directly, which, as a consequence, leads to an algebraic geometry proof of
Theorem~\ref{thm:intro1}. We will pursue this approach in the next section.
Here, we are content to prove the proposition using the symplectic implosion
again, as another example of the power of this construction.

\begin{proof}
  The idea of the proof is to use the fact that $M\times K/K_{\sigma}$ can be
  seen as the symplectic quotient of $M\times \uip$ at any point in $U\subseteq
  \sigma$ with respect to the $T$-action, where $U$ is a neighborhood of
  $\lambda$ in $\sigma$. However, if we recall the definition of $\uip$, we will
  find that a smooth stratum of $M \times \uip$ suffices to give the same
  symplectic quotient. Therefore, instead of using $M \times \uip$ directly, we
  will use a smooth stratum $M \times K/[K_{\sigma},K_{\sigma}] \times \sigma$
  in the following argument without changing the result.

  Let $\Phi_K$ ($\Phi_{T}$ resp.) be the moment map of $M \times
  K/[K_{\sigma},K_{\sigma}] \times \sigma$ of the $K$-action ($T$-action resp.).
  Since $\lambda$ is a regular value of $\Phi_T$, due to \cite[Theorem~7.4 {\&}
  8.10]{Kirwan:1984st}, $T^{\bC}\Phi_T^{-1}(\lambda)$ is the semi-stable point
  set of the $T^{\bC}$-action and
  \begin{equation}
    \label{eq:kt1}
    M \times K/K_{\sigma}\simeq M\times \cO^*_{\lambda} = T^{\bC}\Phi_T^{-1}(\lambda)/T^{\bC}.
  \end{equation}
  Using Lemma~\ref{lm:comm-ss} proved latter, one knows that
  $T^{\bC}\Phi_T^{-1}(\lambda)$ is $G$-invariant, which implies that
  \begin{equation}
    \label{eq:kt2}
    GT^{\bC}\big(\Phi_K^{-1}(0) \cap \Phi_T^{-1}(\lambda)\big)\subseteq GT^{\bC}\Phi_T^{-1}(\lambda) = T^{\bC}\Phi_T^{-1}(\lambda).
  \end{equation}
  By (\ref{eq:kt1}) and (\ref{eq:kt2}), the definition of semi-stable points
  yields
  \begin{equation*}
    (M\times K/K_{\sigma})^{\rss}(\Xi_{\lambda}) = \big(GT^{\bC}(\Phi_K^{-1}(0) \cap \Phi_T^{-1}(\lambda))\big)/ T^{\bC}.
  \end{equation*}
  Similarly, one has
  \begin{equation*}
    (M\times K/K_{\sigma})^{\rss}(\Xi_{\lambda'}) = \big(GT^{\bC}(\Phi_K^{-1}(0) \cap \Phi_T^{-1}(\lambda'))\big)/T^{\bC}.
  \end{equation*}
  
  On the other hand, using the same argument, one can show that
  $G\Phi_K^{-1}(0)$ is $T^{\bC}$-invariant and $\big(GT^{\bC}(\Phi_K^{-1}(0)
  \cap \Phi_T^{-1}(\lambda))\big)/G$ gives the semi-stable point set of the
  $T^{\bC}$-action on $(M \times K/[K_{\sigma},K_{\sigma}] \times \sigma) \dsl
  K$ with respect to the moment map $\mu_{\imp}-\lambda$\footnote{We use the
    same symbol, $\mu_{\imp}$, to denote the moment map of $M_{\imp}$ and its
    restriction on $(M \times K/[K_{\sigma},K_{\sigma}] \times \sigma) \dsl K$}.
  Since $\lambda$ and $\lambda'$ are close enough, by the result of GIT
  quotients for a torus action, the semi-stable point sets of
  $\mu_{\imp}-\lambda$ and $\mu_{\imp}-\lambda'$ coincide,\footnote{Actually,
    this result can be showed by using the flow appearing in the proof of
    Lemma~\ref{lm:ab-eq}.} which implies that
  \begin{equation*}
    GT^{\bC}\big(\Phi_K^{-1}(0) \cap \Phi_T^{-1}(\lambda)\big) = GT^{\bC}\big(\Phi_K^{-1}(0) \cap \Phi_T^{-1}(\lambda')\big).
  \end{equation*}
  As a result, $(M\times K/K_{\sigma})^{\rss}(\Xi_{\lambda}) = (M\times K/K_{\sigma})^{\rss}(\Xi_{\lambda'})$.
\end{proof}

\begin{lemma}
  \label{lm:comm-ss}
  Suppose $K_1$, $K_2$ to be two compact groups. Let $X$ be a Hamiltonian \kah\
  $K_1\times K_2$-manifold (not necessarily compact) and $\mu_1$ be the moment
  map of the $K_1$-action. If $\mu_1^{-1}(0)$ is a compact subset, the set of
  semi-stable point with respect to $\mu_1$, $X^{\rss}(\mu_1)$, is $K_1^{\bC}
  \times K_2^{\bC}$-invariant.
\end{lemma}
\begin{proof}
  By definition, $X^{\rss}(\mu_1)$ is $K_1^{\bC}$-invariant. So we only need to
  show that it is also $K_2^{\bC}$-invariant. Since $\mu_1$ is $K_2$-invariant,
  the definition of semi-stable points implies that $X^{\rss}(\mu_1)$ is
  $K_2$-invariant. Choose a pre-compact open neighborhood $Z\subseteq
  X^{\rss}(\mu_1)$ of $\mu_1^{-1}(0)$. Due to the $K_2$-invariance of
  $X^{\rss}(\mu_1)$, one can find a neighborhood $U$ of the identity element in
  $K_2^{\bC}$ such that $U\cdot Z \subseteq X^{\rss}(\mu_1)$. Using the
  definition of $X^{\rss}(\mu_1)$ again, one has $K_1^{\bC}Z = X^{\rss}(\mu_1)$.
  Therefore,
  \begin{equation*}
    U\cdot X^{\rss}(\mu_1) = UK^{\bC}_1Z = K^{\bC}_1U\cdot Z = X^{\rss}(\mu_1).
  \end{equation*}
  which implies the $K_2^{\bC}$-invariance of
  $X^{\rss}(\mu_1)$ due to the Cartan decomposition $K_2^{\bC} = \exp{\opi
    \kk_2}\cdot K_2$.
\end{proof}

\begin{remark*}
A similar result for the stable point set also holds using a similar argument.
\end{remark*}

\section{The \nth{2} proof of Theorem~\ref{thm:intro1}: the vGIT approach}
\label{sec:vgit}

\subsection{}
In~\cite{Dolgachev:1998wv,Thaddeus:1996tt}, the authors study the variation of
GIT quotients of an algebraic variety when the linearization of the group action
changes. As promised before, in this section, we will use their results to give
another proof of Proposition~\ref{prop:ss} and recover Theorem~\ref{thm:intro1}
consequently. More concretely, we first prove Theorem~\ref{thm:global-prop},
which is merely a restatement of results contained in
Theorem~\ref{thm:dol}~{\&}~\ref{thm:finite}. As a corollary, we can show
Proposition~\ref{prop:ss}. We begin with a discussion about a numerical function
related to the stability condition.

\subsection{A numerical function.}
\label{sub:alg}
We recall some useful definitions of~\cite{Dolgachev:1998wv} in our settings.
Let $(X,\omega_X,\cst_X,\mu_X)$ be a compact Hamiltonian \kah\ $K$-manifold. As
before, we extend the $K$-action on $X$ to a $G$-action holomorphically. A group
homomorphism from $\bC^*$ to $G$ is called a \emph{one-parameter subgroup}, if
it is the complexification of a group homomorphism from $S^1$ to $K$. Naturally,
we can identify such a one-parameter subgroup with an element in $\kk$. Let
$\sX_*(G)\subseteq \kk$ be the set of one-parameter subgroups of $G$. For any
$x\in X$ and $\rho\in \sX_*(G)$, following~\cite{Dolgachev:1998wv}, one defines
a numerical function,
\begin{equation}
  \label{eq:func-m}
  \rM(x) \coloneq \sup_{\rho} d_{\rho}\big(0,\mu_{X}(\overline{\rho(\bC^*)\cdot x})\big),
\end{equation}
in which $d_{\rho}(0,A)$ denotes the signed distance from the origin to the
\underline{\smash{boundary}} of the set $A_{\rho}$. Note that $A_{\rho}\subseteq \bR_+\cdot\rho^*$ is
the projection of $A$ onto the positive ray spanned by $\rho^*$, the dual of $\rho$
under an invariant metric of $\kk$. We remind readers that the definition of
$\rM(x)$ here is a little different from the form given
in~\cite[\S~2.5.2]{Dolgachev:1998wv}. A typo seems to be spotted therein.
Anyway, since our manifold $X$ is only a compact \kah\ manifold, not necessarily
projective, we would like to provide a proof for the following result, which is
well known for the algebraic case.\footnote{After this paper has been submitted,
  we notice the paper~\cite{Georgoulas_2013mo}. And Proposition~\ref{prop:num}
  can be duduced from~\cite[Theorem~14.1]{Georgoulas_2013mo} almostly. We will
  discuss this a little in Remark~\ref{rk:num}.}
\begin{proposition}[A numerical criterion for the stability]
  \label{prop:num}
  Let $X^{\rss}(\mu_X)$ ($X^{\rs}(\mu_X)$ resp.) be the semi-stable (stable resp.)
  point set defined using $\mu_X$ as in \emph{\S}~\ref{sub:ss-point}. With the
  function $\rM(x)$, one can give the following numerical description of
  semi-stable (stable resp.) points, c.f.~\cite[\S~2.5.2]{Dolgachev:1998wv}.
  \begin{align}
    X^{\rss}(\mu_X) &= \{x\in X \,|\, \rM(x) \le 0
                                   \}; \label{eq:m-ss1} \\
    X^{\rs}(\mu_X) &= \{x\in X \,|\, \rM(x) < 0
                                  \}. \label{eq:m-ss2}
  \end{align}
\end{proposition}

\begin{proof}
  Using Atiyah's convexity theorem, which discusses the image of the closure of an
  abelian group action orbit under the moment map,
  c.f.~\cite[Theorem~2]{Atiyah:1982to}, we can reformulate the function $\rM(x)$
  as follows,
  \begin{equation}
    \label{eq:m-1}
    \rM(x) = - \inf_{V\in \sX_*(G)\subseteq \kk} \lim_{t\rightarrow+\infty} \langle\mu_X(\exp{(\opi Vt)}\cdot x),V/\norm{V}\rangle.
  \end{equation}
  The existence of the limit appearing in (\ref{eq:m-1}) is well known,
  c.f.~\cite[\S~3.2]{Mundet-i-Riera:2010aa}. By using (\ref{eq:m-1}), or using
  Atiyah's convexity theorem directly, we can check a special case of
  Proposition~\ref{prop:num}, that is, the case that the group action is
  one-dimensional. In particular, for any $\rho\in \sX_*(G)$,
  Proposition~\ref{prop:num} holds for the $\rho(\bC^*)$-action.

  To prove Proposition~\ref{prop:num} for the general case, we deal with
  (\ref{eq:m-ss1}) and (\ref{eq:m-ss2}) separately. For (\ref{eq:m-ss1}), we can
  use~\cite[Lemma~8.9]{Kirwan:1984st}, which asserts that $x\in X$ is semi-stable
  for the $G$ action if and only if $x$ is semi-stable for every one-parameter
  subgroup of $G$. This result, combined the one-dimensional case that we have
  known, yields (\ref{eq:m-ss1}) immediately.

  To show (\ref{eq:m-ss2}) for general $G$, we need a function
  $\Lambda_x$\footnote{In~\cite[\S~3.3]{Mundet-i-Riera:2010aa}, the same function is
    denoted by $\lambda_x$, we change the notation a little to avoid the symbol
    ambiguity.} used in~\cite{Mundet-i-Riera:2010aa}, also
  c.f.~\cite{Kapovich:2009aa,Woodward:2010aa}. Recall that as a non-positively
  curved space, the symmetric space $\KG$ has a natural compactification by adding
  a \emph{boundary at infinity} $\bKG$. By definition, every point of $z\in \bKG$
  is an equivalent class of geodesics rays on $\KG$. Since the right $G$-action on
  $\KG$ preserves the metric, it induces a right $G$-action on $\bKG$. $\Lambda_x(z)$
  is a Lipschitz continuous function on $\bKG$ with respect to the Tits metric on
  $\bKG$. On $\bKG$, we can define another topology called sphere topology, which,
  in general, is different from the topology induced by the Tits metric. Let $W\in
  \kk$ be a vector of unit norm. The map $\opz:W \mapsto [\exp(\opi Wt)],t\in[0,\infty)$
  leads to a homeomorphism between the unit sphere of $\kk$ and $\bKG$ with sphere
  topology. Using such a homeomorphism, by definition, $\Lambda_x$ can be calculated in
  the following way,
  \begin{equation}
    \label{eq:Lambda}
    \Lambda_x\big(\opz(W)\big)  = \lim_{t\rightarrow+\infty} \langle\mu_X(\exp{(\opi Wt)}\cdot x),W\rangle.
  \end{equation}
  In the following, we need an equivariant property of $\Lambda_x$,
  c.f.~\cite[Lemma~3.6]{Mundet-i-Riera:2010aa}. That is, for any $z\in \bKG$ and
  $g\in G$, we have
  \begin{equation*}
    \Lambda_{g\cdot x}(z) = \Lambda_{x}(z\cdot g).
  \end{equation*}

  We first show the the inclusion in one direction for (\ref{eq:m-ss2}), i.e.
  \begin{equation}
    \label{eq:m-ss21}
    X^{\rs}(\mu_X) \subseteq \{x\in X \,|\, \rM(x) < 0 \}.
  \end{equation}
  Let $x\in X^{\rs}(\mu_X)$. By definition, we can find $y \in \mu_X^{-1}(0)$ such
  that $y = g\cdot x$ and $G_y$ is finite. Then, (\ref{eq:Lambda}) and the
  equivariance of $\Lambda_x$ tell us that for any $V \in \sX_*(G)$,
  \begin{equation*}
    \lim_{t\rightarrow+\infty} \langle\mu_X(\exp{(\opi Vt)}\cdot x),V/ \norm{V}\rangle = \Lambda_x(V/ \norm{V}) = \Lambda_{y} ((V/ \norm{V})\cdot g^{-1}).
  \end{equation*}
  Therefore, by using (\ref{eq:Lambda}) again, to show (\ref{eq:m-ss21}), we
  only need to find $\epsilon_0 > 0$ such that
  \begin{equation*}
    \lim_{t\rightarrow+\infty} \langle\mu_X(\exp{(\opi Wt)}\cdot y),W\rangle \ge \epsilon_0,
  \end{equation*}
  holds for any unit norm vector $W$ in $\kk$. However, by a direct calculation,
  we have the following equality, c.f.~\cite[(3.5)]{Mundet-i-Riera:2010aa},
  \begin{equation*}
    \lim_{t\rightarrow+\infty} \langle\mu_X(\exp{(\opi Wt)}\cdot y),W\rangle = \int^{\infty}_{0} | W^X(\exp{\opi W \tau}\cdot y) |^2 \diff \tau.
  \end{equation*}
  Since $G_y$ is finite, $W^X(y) \ne 0$. Therefore, the existence of $\epsilon_0$ is a
  simple result of the above equality.

  To show the inclusion of (\ref{eq:m-ss2}) in the other direction, we use
  Lemma~\ref{lm:stable-condition}, which is a parallel result of Kirwan's lemma
  for the stable points. As before, such a lemma enables us to reduce the
  general case to the $\bC^*$ case, which we already know. As a result, the proof
  of (\ref{eq:m-ss2}) completes.
\end{proof}

\begin{lemma}
  \label{lm:stable-condition}
  $x \in X$ is a stable point for the $G$-action if and only if for any $\rho\in
  \sX_*(G)$, $x$ is a stable point for the $\rho(\bC^*)$-action with respect to
  restricted moment map.
\end{lemma}
\begin{proof}[Proof of Lemma~\ref{lm:stable-condition}]
  Note that we have proved the inclusion (\ref{eq:m-ss21}). Then the ``only if''
  part the lemma is a consequence of this inclusion and
  Proposition~\ref{prop:num} for the $\bC^*$-action that we have known.

  To show the ``if'' part of the lemma, we first show the following claim:
  \emph{if $x$ is a stable point for any one-parameter subgroup action, then $g\cdot
    x$ satisfies the same property for any $g\in G$.} To prove it, we need
  following result about the group action on $\bKG$,
  c.f.~\cite[\S\S~2,5]{Mundet-i-Riera:2010aa}. For any $V\in \sX_*(G)$, there
  exists $Y\in \sX_*(G)$ such that
  \begin{equation*}
    \opz(Y/\norm{ Y}) = \opz(V/\norm{ V} )\cdot g.
  \end{equation*}
  In other words, the $G$-action on $\bKG$ preserves the ``rational'' points.
  With this property, using the equivariance of $\Lambda_x$ and
  Proposition~\ref{prop:num} for the $\bC^*$ action, the claim follows.

  Now, we can argue by \emph{reductio ad absurdum} to show the ``if'' part of
  the lemma. Assume that $x$ is not stable. By Kirwan's lemma, or
  (\ref{eq:m-ss1}), one knows that $x$ is semi-stable at least. Moreover, since
  $x$ is stable with respect to any one-parameter subgroup action, the isotropic
  subgroup of $x$ is finite. Therefore, supposing $y\in \mu_X^{-1}(0)$ lying in
  the closure of $G\cdot x$, then $y\notin G\cdot x$. On the other hand,
  by~\cite[Corollary~5.5.4]{Woodward:2010aa}, one can find a one-parameter
  group\footnote{In~\cite{Woodward:2010aa}, a one-parameter subgroup does not
    necessarily come from the complexification of an $S^1$ subgroup of $K$. But
    for the one-parameter subgroup $\rho_0$ needed here, by checking the proof
    of~\cite[Corollary~5.5.4]{Woodward:2010aa}, we find that $\rho_0 \in \sX_*(G)$
    indeed.} $\rho_0$ and a point $w\in G\cdot x$ such that $y$ lies in the closure of
  $\rho_0(\bC^*)\cdot w$. Note the claim in the former paragraph implies that $w$ is
  stable with any one-parameter subgroup action due to $w\in G\cdot x$. Especially,
  $w$ is stable with respect to the $\rho_0(\bC^*)$-action. Meanwhile, since $y\in
  \mu_X^{-1}(0)$, it entails that $y$ also lies in the zero level set of the
  moment map associated to the $\rho_0(\bC^*)$-action. Hence, one can conclude that
  $y\in \rho_0(\bC^*)\cdot w \subseteq G\cdot x$, which is a contradiction. Consequently, $x$ must be
  a stable point.
   
\end{proof}

\begin{remark}
  \label{rk:num}
  Some comments about the proof and consequences of Proposition~\ref{prop:num}.
  \begin{enumerate}[label=\Alph*.]
  \item As readers may have noted from the proof of Proposition~\ref{prop:num},
    the stability can be characterized by using the function $\Lambda_x$, which has
    been done in~\cite{Mundet-i-Riera:2010aa}. Similar thing is also true for
    the semi-stability. Using $\Lambda_x$, \cite[Theorem~4.3]{Teleman:2004aa} asserts
    that $x$ is semi-stable if and only if $\Lambda_x(\opz(V))\ge 0$ for any $V\in
    \kk$. Using an argument like~\cite[Remark~5.5.4]{Woodward:2010aa}, this
    result is equivalent to (\ref{eq:m-ss1}). Hence, one can obtain a proof of
    (\ref{eq:m-ss1}) without using Kirwan's lemma.
  \item During the proof of Proposition~\ref{prop:num}, one actually has shown the
    following results,
    \begin{align}
      x\text{ is semi-stable} \Leftrightarrow \lim_{t\rightarrow +\infty}\limits \langle\mu_X(\exp{(\opi Vt)}\cdot x),V\rangle \ge
      0 \text{ for any }V\in \sX_*(G),\label{eq:hm-ss}\\
      x\text{ is stable} \Leftrightarrow \lim_{t\rightarrow +\infty}\limits \langle\mu_X(\exp{(\opi Vt)}\cdot x),V\rangle >
      0 \text{ for any }V\in \sX_*(G),\label{eq:hm-s}
    \end{align}
    which is an analog of the classical Hilbert-Mumford numerical
    criterion~\cite[Theorem~2.1]{Mumford:1994aa}. Such type results have been
    obtained by~\cite{Georgoulas_2013mo}. Moreover, they also prove a similar
    result for the poly-stable points.
  \item If $X$ is projective, Ness has proved a stronger result for the stable
    points set, that is, the infimum in (\ref{eq:m-1}) can be achieved at
    certain $V\in \sX_*(G)$,~\cite[Lemma~3.5]{Ness:1979aa}. But for the \kah\
    case, the same result seems not to be true, even for a torus action.
    However, if we try to take the infimum in $\kk$, it is possible that a
    similar result still holds.
  \item As we have said, in~\cite[Theorem~14.1]{Georgoulas_2013mo}, the authors
    prove the Hilbert-Mumford numerical criterion in a different way. As we have
    seen, up to the infimum in (\ref{eq:m-1}), Proposition~\ref{prop:num} are
    almost equivalent to (\ref{eq:hm-ss}) and (\ref{eq:hm-s}). Their method is
    rather analytic and the key point is an estimate discovered
    by~\cite{Chen_2014ca} using the {\L}ojasiewicz inequality. Compared to their
    method, our method is not so direct but seems to be less technical. But it
    deserves to mention that to show (\ref{eq:hm-s}), both proofs use the
    existence of a slice for the $G$-action in some way.
  \end{enumerate}
\end{remark}

\subsection{A partition of $\sigma$ and the proof of Theorem~\ref{thm:global-prop}.}
\label{sub:wall}
For any $\lambda\in \sigma$, not necessary to be a regular value of the moment map $\mu$ of
$M$, as in \S~\ref{prop:ss}, let $\Xi_{\lambda}$ be the moment map on $M \times K/K_{\sigma}$
induced by the identification between $K/K_{\sigma}$ and $\cO^*_{\lambda}$. To emphasis
its the dependence on $\Xi_{\lambda}$ (or $\lambda$), in the following, we denote the
numerical function $\rM(x)$ on $M \times K/K_{\sigma}$ by $\rM^{\lambda}(x)$. But we should
note that for any given \kah\ form and the corresponding moment map on $M \times
K/K_{\sigma}$, we can always define $\rM^{\bullet}(x)$ for them. That is to say,
$\rM^{\bullet}(x)$ is a function on the space of moment maps on $M \times K/K_{\sigma}$. Via
the moment map $\Xi_{\lambda}$, $\sigma$ can be seen as a convex subset of this space of
moment maps. In this paper, we usually only consider the restriction of
$\rM^{\bullet}(x)$ to $\sigma$.

Following~\cite[\S~3.3]{Dolgachev:1998wv}, we use
$\rM^{\bullet}(x)$ to give a partition of $\sigma\subseteq \kt^*_+$. A subset $H$ of $\sigma$ is
called a \emph{wall} if there exists
\begin{equation*}
  x\in (M\times K/K_{\sigma})_{(>0)} \coloneq \{x|\dim{G_x}> 0\}
\end{equation*}
such that
\begin{equation*}
  H = H(x) \coloneq \{ \lambda\in \sigma| \rM^{\lambda}(x) = 0\}.
\end{equation*}
For a \emph{chamber}, we mean that it is a non-empty connected component of the
complement of the union of walls in $\sigma$. About the relation between
chambers and the sets of semi-stable (or stable) point, one has the following
result. Note that this result, as well as next several results in this section, is
a restatement of the corresponding result in~\cite{Dolgachev:1998wv} for
projective manifolds in our analytic settings. For their proofs, since we have
dropped the algebraicity condition, we will provide details about the necessary
modification compared to their original proofs.

\begin{theorem}[{\cite[Theorem~3.3.2]{Dolgachev:1998wv}}]
  \label{thm:dol}
  Let $\lambda,\lambda'$ be two points in $\sigma$.
  \begin{enumerate}[label=(\emph{\roman*})]
  \item $\lambda$ belongs to some wall if and only if $(M\times
    K/K_{\sigma})^{\rss}(\Xi_{\lambda}) \neq (M\times
    K/K_{\sigma})^{\rs}(\Xi_{\lambda}) $;
  \item $\lambda$ and $\lambda'$ belong to the same chamber if and only if
    \begin{equation*}
      (M\times K/K_{\sigma})^{\rss}(\Xi_{\lambda}) = (M\times
    K/K_{\sigma})^{\rs}(\Xi_{\lambda}) = (M\times
    K/K_{\sigma})^{\rss}(\Xi_{\lambda'}) = (M\times
    K/K_{\sigma})^{\rs}(\Xi_{\lambda'});
    \end{equation*}
  \item each chamber $C$ is convex, and is of the form
    \begin{equation*}
      C = \bigcap_{x\in (M\times K/K_{\sigma})^{\rs}(C)}\{\lambda|\rM^{\lambda}(x)< 0 \},
    \end{equation*}
    where $(M\times K/K_{\sigma})^{\rs}(C) \coloneq (M\times
    K/K_{\sigma})^{\rs}(\Xi_{\lambda})$ for any $\lambda\in C$.
  \end{enumerate}
\end{theorem}

\begin{proof}
  We notice that in our analytic settings, by (\ref{eq:m-1}), $\rM^{\bullet}(x)$,
  viewed as a function on the space of moment maps, is also a sub-additive and
  positively homogeneous function just like its algebraic counterpart
  in~\cite[Lemma~3.2.5]{Dolgachev:1998wv}. It implies that $\rM^{\bullet}(x)$ is a
  convex function. Especially, $\rM^{\bullet}(x)$, viewed as
  function on $\sigma$, is also convex. Given this fact, combined with
  Proposition~\ref{prop:num} and~\cite[Proposition~2.4]{Sjamaar:1995aa}, we can
  repeat all arguments in the proof of~\cite[Theorem~3.3.2]{Dolgachev:1998wv}
  verbatim.
\end{proof}

By Theorem~\ref{thm:dol}, to complete the proof of
Theorem~\ref{thm:global-prop}, therefore the proof of Proposition~\ref{prop:ss},
one needs to show the following finiteness result,
which corresponds to~\cite[Theorem~3.3.3]{Dolgachev:1998wv}.
\begin{theorem}\label{thm:finite}
  There are only finitely many walls in $\sigma$.
\end{theorem}
To prove~\cite[Theorem~3.3.3]{Dolgachev:1998wv}, Dolgachev and Hu use a key
fact,~\cite[Theorem~2.4.5]{Dolgachev:1998wv}\footnote{This result is sometimes
  called Dolgachev-Hu's finiteness theorem in literature, e.g.\
  \cite{Bial-ynicki-Birula:1998aa,Schmitt:2003aa}.}, which asserts that there
are finitely many points, $\lambda_1,\cdots,\lambda_N$ in $\sigma$ such that for
any $\lambda\in \sigma$, the set $(M\times K/K_{\sigma})^{\rss}(\Xi_{\lambda})$
equals to one of the sets $(M\times K/K_{\sigma})^{\rss}(\Xi_{\lambda_i})$. So
if one could show such a result holds for a Hamiltonian \kah\ manifold, one can
repeat Dolgachev and Hu's proof to show Theorem~\ref{thm:finite}. We notice that
although~\cite[Theorem~2.4.5]{Dolgachev:1998wv} is stated not only for the
algebraic case, but also for the general \kah\ case, the proof of which,
however, uses a lemma~\cite[Lemma~1.3.6]{Dolgachev:1998wv}, that is of algebraic
nature.

Here, we try to give a proof of the finiteness theorem of Dolgachev-Hu for \kah\
manifolds by modifying the argument
of~\cite[Example~5.1]{Bial-ynicki-Birula:1998aa}, which gives another proof for
the finiteness theorem in the algebraic case. We begin with a precise statement
of the theorem to be proved.

\begin{theorem}[{\cite[Theorem~2.4.5(ii)]{Dolgachev:1998wv}}]
  \label{thm:finite-ss}
  Suppose that a compact group $K$ acts on a compact complex manifold $X$
  holomorphically. There exist finitely many open subsets of $X$,
  $\{U_1,\cdots,U_M\}$, such that for any Hamiltonian \kah\ structure on $X$
  compatible with the $K$-action, the corresponding semi-stable point set
  must be one of $\{U_1,\cdots,U_M\}$.
\end{theorem}

\begin{proof}
  In spirit, Bia\l{}ynicki's method is similar to the method that we have used
  in the proof of Proposition~\ref{prop:num}, i.e.\ reducing a general compact
  group action to a torus action, and checking the result for the torus action
  using Atiyah's convexity theorem.
  
  \smallskip
  \noindent\textsc{Step 1.} Reducing to the torus action case. As usual, let
  $G=K^{\bC}$ be the complexification of $K$. Choose a maximal torus $T$ in $K$,
  then $T^{\bC}$ is the maximal torus of $G$. For any $K$-invariant Hamiltonian
  \kah\ structure $(\omega,\Psi)$ on $X$, denote the induced moment map for the
  $T$-action by $\Psi_T$. One has the following relation between the semi-stable
  point sets for the $K$ and $T$-actions.
  \begin{equation}
    \label{eq:kt-ss-set}
    X^{\rss}(\Psi) = \bigcap_{k\in K} \big(k\cdot X^{\rss}(\Psi_T)\big).
  \end{equation}
  As in~\cite{Schmitt:2003aa}, (\ref{eq:kt-ss-set}) is a direct result of the
  Hilbert-Mumford criterion, i.e.\ (\ref{eq:hm-ss}), and the equivariance
  property of the moment map. By (\ref{eq:kt-ss-set}), it's clear that if
  Theorem~\ref{thm:finite-ss} holds for the $T$-action, the same theorem also
  holds for the $K$-action.

  \smallskip
  \noindent\textsc{Step 2.} Verifying theorem for the $T$-action. Denote
  $\{F_1,\cdots,F_N\}$ to be the set of connected components of the $T$-action (or
  $T^{\bC}$-action equivalently) fixed points. For any $x\in X$, following
  Bia\l{}ynicki-Birula, we introduce two sets,
  \begin{equation}
    \label{eq:sh-c}
    \opsh{(x)} \coloneq \{F_i | \overline{T^{\bC}\cdot x} \cap F_i \neq \emptyset\},\; \opc{(x)} \coloneq \{y\in X |
    \opsh{(x)} = \opsh{(y)} \}.
  \end{equation}
  By definition, there are only finitely many different sets lying in
  $\{\opc{(x)} | x\in X \}$. And for any $x,y\in X$, one has either $\opc{(x)} =
  \opc{(y)}$ or $\opc{(x)} \cap \opc{(y)} = \emptyset$.
  
  For any \kah\ form $\omega_X$ and moment map $\mu_X$ on $X$,
  by~\cite[Theorem~2]{Atiyah:1982to}, i.e.\ Atiyah's convexity theorem, the
  semi-stable point set has the following representation,
  \begin{equation}
    \label{eq:ss-t}
    X^{\rss}(\mu_X) = \{x\in X | 0\in \mu_X(\overline{T^{\bC}\cdot x})\} = \{x\in X | 0\in
    \conv{\big(\mu_X(\opsh (x))\big)} \}.
  \end{equation}
  By (\ref{eq:sh-c}) and (\ref{eq:ss-t}), there exist $x_1,\cdots,x_k\in X$ such that
  $X^{\rss}(\mu_X) = \sqcup_{i=1}^k \opc{(x_i)}$, which implies that there are only
  finite possibilities for the shape of $X^{\rss}(\mu_X)$.

\end{proof}

\begin{remark}
  \label{rk:finite-ss}
  Some comments about Theorem~\ref{thm:finite-ss}.
  \begin{enumerate}[label=\Alph*.]
  \item For the algebraic case, the theorem holds not only for smooth manifolds
    but also for projective varieties, c.f.~\cite{Bial-ynicki-Birula:1998aa,
      Dolgachev:1998wv}. In fact, the theorem is even true for the positive
    characteristic, c.f.~\cite{Ressayre:2000aa,Schmitt:2003aa}. Taking this fact
    into consideration, it seems reasonable to expect a similar result to hold
    for compact Hamiltonian~\kah\ spaces. As before, to obtain such a
    generalization, we need two main results for singular spaces. One of them,
    i.e. Atiyah's convexity theorem, does have a singular space
    generalization,~\cite[p.80,~Theorem]{Heinzner:1996aa}. So, the only thing
    left is to find a \kah\ spaces version of the Hilbert-Mumford criterion.
    
  \item We notice that an argument similar to that of
    Theorem~\ref{thm:finite-ss} also appears in~\cite[pp.~174,
    Proposition]{Heinzner:2001aa}. Moreover, by the result of that paper, on a
    projective manifold, any open set $U_i$ appearing in
    Theorem~\ref{thm:finite-ss} comes actually from the corresponding set
    constructed using the algebraic method.

  \item In~\cite[Theorem~3.5]{Hu_1995aa}, the author uses a method very close to
    the method used in the \textsc{Step 2.} in the above proof. In particular,
    this means that although~\cite{Hu_1995aa} discusses only about projective
    manifold, this beautiful result is also true for general compact \kah\ manifold.
  \end{enumerate}
\end{remark}

\subsection{}
With Theorem~\ref{thm:finite-ss} in hand, one can use the same argument as
in~\cite[Theorem~3.3.3]{Dolgachev:1998wv} to prove Theorem~\ref{thm:finite}. For
readers' convenience, we incorporate Dolgachev and Hu's proof here.
\begin{proof}[Proof of Theorem~\ref{thm:finite}]
  For any wall $H$, let
  \begin{equation*}
    (M\times K/K_{\sigma})_H \coloneq \{x\in (M\times K/K_{\sigma})_{(>0)} | H \subseteq H(x)\},
  \end{equation*}
  where $(M\times K/K_{\sigma})_{(>0)}$ is the set of points in $M\times
  K/K_{\sigma}$ with a positive-dimensional isotropic group and $H(x)$ is the set
  defined in \S~\ref{sub:wall}. We will show that $(M\times K/K_{\sigma})_H$ can
  be written in the following way,
  \begin{equation}
    \label{eq:wall}
    (M\times K/K_{\sigma})_H = \bigcap_{\lambda\in H}\big((M\times K/K_{\sigma})^{\rss}
    (\Xi_{\lambda}) \cap (M\times K/K_{\sigma})_{(>0)}\big).
  \end{equation}

  On one side, if $x\in (M\times K/K_{\sigma})_H$, one has $H\subseteq H(x)$ by
  definition, which implies that, due to the definition of $H(x)$,
  $\rM^{\bullet}(x)$ vanishes on $H$. Then by Proposition~\ref{prop:num}, for
  any $\lambda\in H$, $x\in (M\times K/K_{\sigma})^{\rss}(\Xi_{\lambda})$, i.e.\
  \begin{equation*}
    x\in \bigcap_{\lambda\in H} \big((M\times K/K_{\sigma})^{\rss}(\Xi_{\lambda})
    \cap (M\times K/K_{\sigma})_{(>0)}\big).
  \end{equation*}
  On the other side, if
  \begin{equation*}
    x\in \bigcap_{\lambda\in H}\big((M\times K/K_{\sigma})^{\rss}(\Xi_{\lambda}) \cap
    (M\times K/K_{\sigma})_{(>0)}\big),
  \end{equation*}
  by Proposition~\ref{prop:num} again, for
  any $\lambda\in H$, one has $\rM^{\lambda}(x) = 0$, or equivalently
  $\lambda\in H(x)$. Hence, $H\subseteq H(x)$, i.e. $x\in (M\times
  K/K_{\sigma})_H$. As a result, (\ref{eq:wall}) is true.

  By Theorem~\ref{thm:finite-ss}, one can find a finite set of points
  $\{\lambda_1,\cdots, \lambda_M\}\subseteq \sigma$ such that for any
  $\lambda\in\sigma$, the set $(M\times K/K_{\sigma})^{\rss}(\Xi_{\lambda})$
  equals to one of the sets $(M\times K/K_{\sigma})^{\rss}(\Xi_{\lambda_i})$. By
  (\ref{eq:wall}), we know that there are only finitely many subsets of $M\times
  K/K_{\sigma}$ which are of the form $(M\times K/K_{\sigma})_H$ for some wall
  $H$. However, by the definition of walls, two walls $H,H'$ coincide if and
  only if $(M\times K/K_{\sigma})_H = (M\times K/K_{\sigma})_{H'}$. As a
  consequence, only finitely many walls exist.
\end{proof}

\section{A relation between two approaches}
\label{sec:comp-betw-two}

\subsection{}
After giving two proofs of Theorem~\ref{thm:intro1}, one by the symplectic
implosion, one by the vGIT, it is appropriate to have a comparison between two
approaches. As we have declared in Introduction, the general guideline behind
two approaches is a reflection of the same plain idea: reducing a non-abelian
reduction problem to an abelian one, which looms in the statement and proof of
Theorem~\ref{thm:introeq} and Theorem~\ref{thm:finite-ss} for example. Besides,
one can go beyond such a general discussion and work out a more concrete
relation between these two approaches. As it turns out, to carry out such a
comparison, it is Thaddeus's proof of vGIT~\cite{Thaddeus:1996tt} that fits better
for this purpose, although we have used Dolgachev and Hu's proof of the same
theory extensively in the previous section. More precisely, an interesting
construction used by Thaddeus~\cite[\S~3.1]{Thaddeus:1996tt} has a natural
correspondence with the symplectic implosion construction. We summarize the
results in Proposition~\ref{prop:ms-imp}~\&~\ref{prop:master-space}.

\subsection{Thaddeus's master space.}
\label{sub:def-ums}
As before, we assume that $(X,\omega_X,\cst_X,\mu_X)$ is a compact Hamiltonian
\kah\ manifold. But here, we further require that there is a holomorphic
prequantum line bundle $L$ over $X$. In other words, the $K$-action is lifted to
a group action on $L$ and there is a $K$-invariant Hermitian metric $h$ on $L$.
Denoting the Chern connection of $h$ by $\nabla$, one requires that
$\frac{\nabla^2}{2\pi \opi} = \omega_X$. Since the $K$-action on $L$ can be
extended to a $G$-action on $L$, in the GIT terminology, $L$ is a
$G$-linearization of $X$.

Recall that $\Pi = \{\varpi_1,\cdots, \varpi_r \}$ is the set of fundamental
weights of $K$. And we denote $\bC_i$, $i=1,\cdots,r$, to be the 1-dimensional
representation of $T$ with weight $\varpi_i$. Also, recall that $V_i$ is the
irreducible representation of $K$ with the highest weight $\varpi_i$ and $v_i$
is a highest weight vector of $V_i$. Then we can and will identify $\bC_i$ with
$\bC v_i\subseteq V_i$. Moreover, let $\varpi_0$ be the zero weight and $\bC_0$
be the trivial representation of $T$. As usual, by using the Borel subgroup
$T^{\bC}\subseteq B \subseteq G$ corresponding to the positive roots we have
chosen, one constructs a line bundle $L_i = G \times_B \bC_i$, $i=0,\cdots,r$,
over $G/B$, associated to the $T^{\bC}$ (or $B$) representation $\bC_i$. By
choosing a parameter $t\in \Delta= \{(t_i) \in \bQ^{r+1}_{\ge 0} | \sum_{i=0}^r
t_i = 1\}$, one has a family of $G$-linearizations, $L\otimes \prod_{i=0}^r
L_i^{-t_i}$\footnote{Without ambiguity, we don't distinguish a bundle on $X$ or
  $G/B$ and its pull-back on $X \times G/B$.}, over $X\times G/B$.

Roughly speaking, Thaddeus~\cite{Thaddeus:1996tt} constructs a ``master space'' $X_{\ms}$, that
transforms GIT quotients of $X\times G/B$ defined by the $G$-action on $L\otimes \prod_{i=0}^r
L_i^{-t_i}$ to GIT quotients of $X_{\ms}$ defined by a family of $T^{\bC}$-actions
on a fixed line bundle. In our situations, one can even construct certain
universal master space, that is, a space independent of $X$, as follows. For
$E^N = \oplus_{i=1}^r \bC v_i \subseteq E = \oplus_{i=1}^r V_i$,
\begin{equation*}
  \ums \coloneq G \times_B \bP(\bC_0 \oplus E^N) = K \times_T \bP(\bC_0 \oplus E^N).
\end{equation*}
The associated (relative) hyperplane line bundle of $\ums$ is denoted
by $\cO_{\ums}(1)$. $\ums$ has the following embedding,
\begin{equation*}
\begin{array}{cccc}
  (\opp,\pi):&\ums & \rightarrow & \bP(\bC_0 \oplus E) \oplus G/B\\
  & [g,[u]] & \mapsto & ([g\cdot u], [g]).
\end{array}
\end{equation*}
where $g\in G,u\in \bC_0 \oplus E^N \subseteq \bC_0 \oplus E$. One has
$\cO_{\ums}(1) = \opp^*\cO_{\bP(\bC_0 \oplus E)}(1)$.

\uline{A notation warning}: $\bP$ used in this paper is the projectivization of
a space, that is, taking all one-dimensional subspaces of the original space,
which is different from the algebraic usage of $\bP$ as
in~\cite{Thaddeus:1996tt}. In fact, the description of the construction given
here is equivalent to Thaddeus's original construction except that we choose to
perform the manipulation on the \uline{dual} bundle. As a useful notation, we
also denote $\bP(\bC_0 \oplus E)$ ($\bP(\bC_0 \oplus E^N)$ resp.) by $\bar{E}$
($\bar{E}^N$ resp.) in the following.

Up to now, we still miss an important assumption in Thaddeus's construction, that
is, the line bundles over $X \times G/B$ used in the construction should be
ample. In our description of the construction, it's equivalent to the negativity
of the line bundles $\{L_i\}$. However, by some direct calculations, one knows
that $L_i$ is only semi-negative for any $i\in\{0,\cdots,r\}$. One can remedy
this problem as follows. Let $\epsilon$ be a small fractional weight lying in
$(\kt^*_+)^{\circ}$. Instead of using a set of weights
$\{\varpi_0,\varpi_1,\cdots,\varpi_r\}$ to construct $\ums$, we use $
\{\varpi_0+\epsilon,\varpi_1+\epsilon,\cdots,\varpi_r+\epsilon\}$ to do the real
job.\footnote{Admittedly, there are no representations associated to such
  fractional weights. One can either use scaling arguments to correct this
  little flaws or take the following discussion of $\ums^{\epsilon}$ and
  $\cO_{\ums^{\epsilon}}(1)$ as definitions directly.} The resulting space is
denoted by $\ums^{\epsilon}$. It turns out that $\ums^{\epsilon}$ is
holomorphically isomorphic to $\ums$. But the (relative) hyperplane line bundle
is changed to (an ample line bundle) $\cO_{\ums^{\epsilon}}(1) = \cO_{\ums}(1)
\otimes \pi^*L^{-1}_{\epsilon}$, where $L_{\epsilon}$ is the line bundle over
$G/B$ associated to $\epsilon$. The first Chern form of
$\cO_{\ums^{\epsilon}}(1)$, i.e.\ the symplectic form on $\ums^{\epsilon}$, is
$\opp^*\omega_{\bar{E}} + \pi^*\cO^*_{\epsilon}$, where $\omega_{\bar{E}}$ is
the Fubini-Study form on $\bar{E}$ and we use the identification $G/B \simeq K/T
\simeq \cO^*_{\epsilon}$. The \emph{master space} of $X$\footnote{If we stick to
  Thaddeus's terminology, it may be more proper to call $(X \times
  \ums^{\epsilon})\dsl K$ the master space of $X \times \ums^{\epsilon}$.} is
defined to be $X_{\ms}\coloneq (X \times \ums^{\epsilon})\dsl K$.

\subsection{}
By the definition of $\ums^{\epsilon}$, there exists a natural $G$ (or $K$)
-action on it and this action can be lifted to an action on
$\cO_{\ums^{\epsilon}}(1)$. About the $T^{\bC}$ (or $T$) -action on
$\ums^{\epsilon}$, one can proceed resembling the symplectic implosion
construction. Namely, we define a $T^{\bC}$-action on $E^N$ by requiring that
the $T^{\bC}$-action is diagonalized and the weight of the action on $\bC_i$ is
given by $-\varpi_i$, $i=0,1,\cdots,r$, which induces the $T^{\bC}$-action on
$\ums^{\epsilon} = \ums$ and $\cO_{\ums}(1)$. As for the $T^{\bC}$-action on
$\cO_{\ums^{\epsilon}}(1)$, one should also take the natural $T^{\bC}$-action on
$L_{\epsilon}$ into consideration. Denote the moment map of the $T$-action on
$\cO_{\ums^{\epsilon}}(1)$ by $\vartheta_T$.

For concrete calculations, it is convenient to identify $T^{\bC}$ with
$(\bC^*)^r$ by using the chosen set of fundamental weights $\Pi$. Especially,
under such an identification, for $z = (z_1,\cdots,z_r)\in (\bC^*)^r\simeq
T^{\bC}$ and $u = (u_1,\cdots,u_r)\in E^N$, one has $z\cdot u =
(z_1^{-1}u_1,\cdots,z_r^{-1}u_r)$.

\subsection{}
\label{sub:compare}
In~\cite[Proposition~7.7]{Guillemin:2002sx}, the authors construct a smooth
manifold closely related to $\ums^{\epsilon}$, i.e.\ $G\times_B E^N$ (the
symplectic form on it depends on $\epsilon$), as a desingularization of the
universal implosion section $\uip = G_N$. Let $\inc:E^N \rightarrow \bar{E}^N$
be the inclusion map, $\inc(u)\coloneq [1:u]$, which induces a $K\times
T$-equivariant inclusion from $G\times_B E^N$ to $\ums^{\epsilon}$, also denoted
by $\inc$. Therefore, $\ums^{\epsilon}$ is a smooth compactification of
$G\times_B E^N$. On the other hand, one notices that the map $\inc$ is the
restriction of the map $\Inc: E \rightarrow \bar{E}$. Let $\overline{G_N}$ be
the closure of $\Inc(G_N)$ in $\bar{E}$ (with respect to either standard or
Zariski topology). $\ums^{\epsilon}$ can also be seen as a desingularization of
$\overline{G_N}$. In fact, the following diagram commutes, where $\opq:G\times_B
E^N \rightarrow E$ is the multiplication.
\begin{equation*}
  \begin{tikzcd}
    G\times_B E^N \arrow[r, hook, "\inc"] \arrow[d, two heads, "\opq"'] & \ums^{\epsilon}
    \arrow[d, two heads, "\opp"] \\
    G_N \arrow[r,hook, "\Inc"] & \overline{G_N}
  \end{tikzcd}.
\end{equation*}

Now by using some scaling, we assume that $\mu_X(X) \cap \kt^*_+$ is contained
in the simplex spanned by vertices $\{0,\varpi_1/\sqrt{2\pi r},\allowbreak
\cdots, \varpi_r/\sqrt{2\pi r}\}$ and $\epsilon\in (\kt^*_+)^{\circ} \cap
\mu_X(X)$ is a small regular value of $\mu_X$. One can check that, under such an
assumption, $X_{\ms}$ is an orbifold. Meanwhile, due
to~\cite[Corollary~7.13]{Guillemin:2002sx}, $(X \times G\times_B E^N)\dsl K$ is
also an orbifold, which in fact is a partial desingularization of $X_{\imp}$,
that is, there exists a proper surjective bimeromorphic map $(X \times G\times_B
E^N)\dsl K \rightarrow X_{\imp}$. Hence, the following proposition unveils the
close relation between the master space and the symplectic implosion, i.e.
$X_{\ms}$ is a partial desingularization of $X_{\imp}$.
\begin{proposition}
  \label{prop:ms-imp}
  $(X \times G\times_B E^N)\dsl K$ is holomorphically isomorphic to $X_{\ms}$.
\end{proposition}

\begin{proof}
  We denote the moment maps on $X\times G\times_B E^N$ and $X \times
  \ums^{\epsilon}$ by $\Phi,\Psi$. We are going to show that $(X\times G\times_B
  E^N)^{\rss}(\Phi)$ is isomorphic to $(X\times \ums^{\epsilon})^{\rss}(\Psi)$
  under the map $\id_X \times \inc$.

  Let $(x,[e,u])\in X\times G\times_B E^N$ satisfy $\Phi(x,[e,u]) = 0$. We claim
  that there exists $z \in (\bC^{r})^*\simeq T^{\bC}$ such that
  $\Psi(x,[e,\inc{(z\cdot u)}]) = 0$. By (\ref{eq:vectsymp})
  and~\cite[\S~2.7]{Kirwan:1984st}, for $u =(u_1,\cdots,u_r) \in
  \oplus_{i=1}^r\bC_i = E^N \subseteq E$, one can calculate the moment maps of
  $E$ and $\bar{E}$ as follows,\footnote{A minus sign appears in
    $\mu_{\bar{E},K}$ compared to Kirwan's formula, since Kirwan's sign
    convention of a moment map is different from us.}
  \begin{equation}
    \label{eq:two-mm}
    \begin{aligned}
      \mu_{E,K}(u) &= -\pi \sum_{i=1}^r |u_i|^2\varpi_i,\\
      \mu_{\bar{E},K}\big(\inc{(u)}\big) &= \frac{-1}{1 + \norm{u}^2} \sum_{i=1}^r |u_i|^2\varpi_i.
    \end{aligned}
  \end{equation}
  Choose
  \begin{equation*}
    z = \sqrt{\frac{1- \pi \sum^{r}_{i=1}|u_i|^2}{\pi}}(1,\cdots,1),
  \end{equation*}
  which is well-defined due to the assumption on the image of $\mu_X$. By
  (\ref{eq:two-mm}), one can check that $\mu_{E,K}(u) =
  \mu_{\bar{E},K}\big(\inc{(z\cdot u)}\big)$. Therefore, by the definition of
  the moment maps of $G\times_B E^N$ and $\ums^{\epsilon}$ in
  \S~\ref{sub:def-ums} and~\cite[Proposition~7.7]{Guillemin:2002sx}, one has
  \begin{equation*}
    \Psi(x,[e,\inc{(z\cdot u)}]) = \mu_X(x) + \mu_{\bar{E},K}\big(\inc{(z\cdot u)}\big) -\epsilon = \mu_X(x) +
    \mu_{E,K}(u) - \epsilon = \Phi(x,[e,u]) = 0.
  \end{equation*}
  Hence,
  \begin{equation*}
    z\cdot\big(x,\inc{([e,u])}\big) = (x,[e,\inc{(z\cdot u)}]) \in
    (X\times \ums^{\epsilon})^{\rss}(\Psi),
  \end{equation*}
  which implies
  $\big(x,\inc{([e,u])}\big) \in (X\times \ums^{\epsilon})^{\rss}(\Psi)$ due to
  Lemma~\ref{lm:comm-ss}. As a result, using the equivarince property of the map
  $\inc$, we have the inclusion for one direction:
  \begin{equation*}
    \id \times \inc\big((X\times G\times_B E^N)^{\rss}(\Phi)\big) \subseteq (X\times
    \ums^{\epsilon})^{\rss}(\Psi).
  \end{equation*}

  To show the inclusion of the other direction, one only needs to notice the
  following fact: for any point $p\in \Psi^{-1}(0)$, there exists $k\in K$ such
  that $k\cdot p$ is of the form $(x,[e,[1:u]])$ and lying in the image of $\id
  \times \inc$ consequently. To verify this fact, we note that for any point
  $p\in \Psi^{-1}(0)$ that does not have the asserted property, $p$ must be of the
  form $(x,[e,[0:u]])$ up to a $K$-action on $p$. But due to the assumption on
  the image of $\mu_X$ again,
  \begin{equation*}
    \Psi(x,[e,[0:u]]) = \mu_X(x) + \epsilon -
    \frac{1}{\norm{u}^2} \sum_{i=1}^r |u_i|^2\varpi_i
  \end{equation*}
  can never vanish, which is
  a contradiction. With such an result, we are now in essentially the same
  situation as in the previous paragraph, which enables us to use the same
  argument to show that $\id \times \inc\big((X\times G\times_B
  E^N)^{\rss}(\Phi)\big) \supseteq (X\times \ums^{\epsilon})^{\rss}(\Psi)$.
\end{proof}
\begin{remark*}
  In a sense, the above proposition compares two kinds of symplectic quotients,
  coming from an affine manifold and its embedding into a projective space
  respectively. One notices that similar results have been discussed
  in~\cite[pp.~224, 242]{Kirwan:2011aa}.
\end{remark*}

\subsection{Torus actions on the master space.}
\label{sub:sit-action}
In~\cite{Thaddeus:1996tt}, the author introduces a family of torus actions on
$\cO_{\ums^{\epsilon}}(1)$ and calculates the corresponding GIT quotients. Not
unexpected, such quotients can also be obtained as symplectic quotients with
respect to different values of a moment map. To make this comparison more
transparent, we first recall the definition of the family of torus actions.

Let $\bfT = \{\xi\in \bC^{r+1}| \prod_{i=0}^r \xi_i =1\}$. Then $\bfT$ acts on
\begin{equation*}
  (u_0^{\epsilon},\cdots, u_r^{\epsilon})\in (\bC_0 \oplus E^N)\otimes
  \bC_{\epsilon} = \oplus_{i=0}^r \bC_i\otimes \bC_{\epsilon}
\end{equation*}
by
\begin{equation*}
  \xi\cdot(u_0^{\epsilon},\cdots,u_r^{\epsilon}) \coloneq
  (\xi_0u_0^{\epsilon},\cdots ,\xi_ru_r^{\epsilon}),
\end{equation*}
which induces a
$\bfT$-action on $\ums^{\epsilon}$ and $\cO_{\ums^{\epsilon}}(1)$. Now, recall
that one has a parameter $t\in \Delta$, which determines a fractional character
of $\bfT$ by $t(\xi) \coloneq \prod_{i=0}^r \xi_i^{t_i}$. Using this, the
claimed family of $\bfT(t)$-actions on $\cO_{\ums^{\epsilon}}(1)$ is induced
from the following family of actions depending on $t$:
\begin{equation*}
  \xi\cdot_t (u_0^{\epsilon},\cdots,u_r^{\epsilon}) \coloneq
  t^{-1}(\xi)(\xi_0u_0^{\epsilon},\cdots ,\xi_ru_r^{\epsilon}).
\end{equation*}

Using the identification between $T^{\bC}$ and $(\bC^*)^r$, one defines a map from
$T^{\bC}$ to $\bfT$ as follows. Let $\eta = (\prod_{i=1}^rz_i)^{\frac{1}{r+1}}$.
\begin{equation*}
\begin{array}{cccc}
    \varphi:& T^{\bC}\simeq (\bC^*)^r & \rightarrow & \bfT\\
     & (z_1,\cdots,z_r) &\mapsto & (\eta, \eta z_1^{-1},\cdots,\eta z_r^{-1}).
\end{array}
\end{equation*}
Strictly speaking, $\varphi$ is not a map, which is just a tuple of fractional
characters of $(\bC^*)^r$. But for our following usage, such a ``map'' $\varphi$
is enough. Therefore, we will ignore this little inaccuracy in the definition of
$\varphi$.

By composing $\varphi$ and the $\bfT$-action, one has a new $T^{\bC}$-action on
$\cO_{\ums^{\epsilon}}(1)$, whose $T$-moment map is
\begin{equation*}
  \vartheta'_T = \vartheta_T - \frac{1}{r+1}\sum_{i=0}^{r}\varpi_i -\epsilon.
\end{equation*}
In the same way, by composing
$\varphi$ and the $\bfT(t)$-action, one has a family of $T^{\bC}(t)$-actions
depending on $t\in\Delta$ on $\cO_{\ums^{\epsilon}}(1)$, whose $T(t)$-moment
maps are
\begin{equation*}
  \vartheta'_T(t) = \vartheta_T - \sum_{i=0}^{r}t_i(\varpi_i+\epsilon).
\end{equation*}
Therefore, using the Kirwan-Ness theorem, we come to the following
interpretation of the torus action quotients on a master space.

\begin{proposition}
  \label{prop:master-space}
  Under the map $\varphi$, for any $t\in \Delta$, the GIT quotient $X_{\ms}\dsl \bfT(t)$ is
  holomorphically isomorphic to the symplectic quotient
  $X_{\ms}\dsl_{\sum_{i=0}^{r}t_i(\varpi_i+\epsilon)} T$.
\end{proposition}

In fact, one can also calculate and compare the two quotients explicitly as
follows. On the algebraic side, by results in~\cite[\S~3.1]{Thaddeus:1996tt},
$X_{\ms}\dsl \bfT(t)$ is isomorphic to $(X\times G/B)\dsl G(t)$, which is the
GIT quotient with respect to the $G$-linearization $L\otimes \prod_{i=0}^r (L_i
\otimes L_{\epsilon})^{-t_i}$. On the symplectic side, by
Proposition~\ref{prop:ms-imp},
\begin{equation*}
  X_{\ms}\dsl_{\sum_{i=0}^{r} t_i(\varpi_i+\epsilon)} T = \big((X \times G\times_B E^N)\dsl K\big)
\dsl_{\sum_{i=0}^{r}t_i(\varpi_i+\epsilon)} T = (X \times \cO^*_{\sum_{i=0}^{r}
  t_i(\varpi_i+\epsilon)}) \dsl K,
\end{equation*}
which is isomorphic to $(X\times G/B)\dsl
G(t)$ by the Kirwan-Ness theorem again.

Proposition~\ref{prop:ms-imp} and Proposition~\ref{prop:master-space} complete
our interpretation of Thaddeus's construction in terms of the symplectic
implosion techniques.

\bibliographystyle{amsplain}

\end{document}